 \documentclass[a4paper]{amsart}

\usepackage[utf8]{inputenc}
\usepackage{t1enc}

\usepackage{amssymb}
\usepackage[all]{xy}

\usepackage{mathbbol}
\usepackage{hyperref}

\usepackage{latexsym}
\usepackage[mathscr]{euscript}
\newcommand{\comments}[1]{}
\newcommand{\cst}{\textup{C}^*}

\newcommand{\id}[1]{\operatorname{id}_{#1}}
\newcommand{\I}[1]{\operatorname{I}_{#1}}
\newcommand{\Mor}{\operatorname{Mor}}

\newcommand{\Rep}{\operatorname{Rep}}

\newcommand{\flip}{\operatorname{flip}}

\newcommand{\M}{\operatorname{M}}

\newcommand{\B}{\operatorname{B}}
\newcommand{\tens}{\otimes }

\newcommand{\Ahat}{{\widehat A}}
\newcommand{\omegahat}{{\widehat\omega}}

\newcommand{\ahat}{{\widehat a}}

\newcommand{\Mhat}{{\widehat M}}

\newcommand{\Ghat}{{\widehat G}}

\newcommand{\Rhat}{{\widehat R}}
\newcommand{\Vhat}{{\widehat V}}
\newcommand{\Deltahat}{{\widehat\Delta}}

\newcommand{\rhohat}{{\widehat\rho}}
\newcommand{\sigmahat}{{\widehat\sigma}}

\newcommand{\alphahat}{{\widehat\alpha}}
\newcommand{\pihat}{\widehat{\pi}}

\newcommand{\ehat}{\widehat{e}}

\newcommand{\xbar}{\overline{x}}

\newcommand{\Hbar}{\overline{H}}

\newcommand{\varphitilde}{{\widetilde\varphi}}

\newcommand{\alphatilde}{{\widetilde \alpha}}

\newcommand{\pitilde}{{\widetilde \pi}}

\newcommand{\psitilde}{\widetilde{\psi}}

\newcommand{\natu}{{\mathbb N}}

\newcommand{\inte}{{\mathbb Z}}

\newcommand{\Vs}[1]{\rule{0mm}{#1mm}}
\newcommand{\vs}{\vspace{2mm}}

\newcommand{\compcent}[1]{\vcenter{\hbox{$#1\circ$}}} 
\newcommand{\comp}{\mathbin{\mathchoice 
{\compcent\scriptstyle}{\compcent\scriptstyle} 
{\compcent\scriptscriptstyle}{\compcent\scriptscriptstyle}}}

\newcommand{\rf}[1]{{\rm (\ref{#1})}}

\newcommand{\set}[2]{\left\{#1:#2\right\}}
\newcommand{\Vset}[2]{\left[#1:#2\right]}

\newcommand{\norm}[1]{\left\|#1\right\|}

\newcommand{\ITS}[3]{\left(#1\,\vline\,#2\,\vline\,#3\right)}
\newcommand{\bra}[1]{\left(#1\,\right\vert}
\newcommand{\ket}[1]{\left\vert#1\right)}
\newcommand{\xxx}[1]{}

\numberwithin{equation}{section}
\newtheorem{Thm}{Theorem}[section]
\newtheorem{Def}[Thm]{Definition}

\newtheorem{Prop}[Thm]{Proposition}

\newtheorem{rem}[Thm]{Remark}

\newenvironment{pf}[1][]{\vs
\noindent{\it Proof#1. }}{\qed}

\newcommand{\etyk}[1]{\label{#1}\stepcounter{equation}\tag{\theequation}}

\usepackage{xcolor}
\newcommand{\SLW}[1]{{\color{red}#1}}

\begin{document}

\title[Landstad-Vaes theory]{Landstad-Vaes theory for \\locally compact quantum groups}
\author{Sutanu Roy}
\email{sutanu@niser.ac.in}
\address{School of Mathematical Sciences,\\ National Institute of Science Education and Research\\ Bhubanes\-war, HBNI, Jatni, India}

\author{Stanisław Lech Woronowicz}
\email{Stanislaw.Woronowicz@fuw.edu.pl}
\address{Institute of Mathematics of the Polish Academy of Sciences, ul.~Śniadeckich 8, 00-656 Warszawa, Poland and Department of Mathematical Methods in Physics, Faculty of Physics, University of Warsaw, ul.~Pasteura 5, 02-093 Warszawa Poland.}

\date{}

\begin{abstract}
Landstad-Vaes theory deals with the structure of the crossed product of a $\cst$-algebra by an action of locally compact (quantum) group. In particular it describes the position of original algebra inside crossed product. The problem was solved in 1979 by Landstad for locally compact groups and in 2005 by Vaes for regular locally compact quantum groups. To extend the result to non-regular groups we modify the notion of $G$-dynamical system introducing the concept of weak action of quantum groups on $\cst$-algebras. It is still possible to define crossed product (by weak action) and characterise the position of original algebra inside the crossed product. The crossed product is unique up to an isomorphism. At the end we discuss a few applications. 
\end{abstract}


\thanks{Supported by the Alexander von Humboldt-Stiftung. Sutanu Roy was partially supported by Inspire faculty award given by D.S.T., Government of India. S.L. Woronowicz was partially supported by the National Science Center (NCN) grant no.~2015/17/B/ST1/00085.}

\maketitle

\section{Introduction}
\label{sek0}

The concept of crossed product of a $\cst$-algebra by an action of a locally compact (quantum) group comes from the desire to unite in a single object the $\cst$-algebra and the unitaries implementing the action. For the theory of crossed product of  $\cst$-algebras by actions of locally compact groups see \cite{Ped}.\vs

One of the most used formula in the operator algebra theory is the implementing of an automorphism by a unitary operator:
\[
\etyk{imp}
\alpha(d)=UdU^{*}.
\]
In this formula $U$ is a unitary operator acting on a Hilbert space $K$ and $d$ runs over a $\cst$-algebra $D$ of operators acting on $K$. It is assumed that $\alpha(d)\in D$ for any $d\in D$ and that $\alpha(D)=D$. Then $\alpha$ is an automorphism of $\cst$-algebra $D$. We say that automorphism $\alpha$ is implemented by $U$.\vs

Any automorphism  of a $\cst$-algebra $D$ can be implemented. For any $\alpha$ there exist a Hilbert space $K$ and a pair $(j,U)$, where $j$ is an embedding of $D$ in $\B(K)$ and $U$ is a unitary operator acting on $K$ such that (identifying $d$ with $j(d)$) we have \rf{imp}. We say that $(j,U)$ is a covariant representation of $D$.\vs

It is interesting to extend $D$ by including $U$. Let
\[
\etyk{1744}
B=\set{\Vs{3.5}U^{n}d}{n\in\inte,d\in D}^{\rm CLS},
\]
where ${\rm CLS}$ stays for norm closed linear span. Then $B$ is a $\cst$-algebra, $U\in\M(B)$, $D\subset B$ and $DB=B$. The latter means that the embedding $D\subset B$ is a morphism from $D$ into $B$. %We refer to the next section for the concept of morphism in the category of $\cst$-algebras. 
In general, for given $D$ and $\alpha$, the algebra $B$ may depend on the used covariant representation. For instance, if $\alpha$ is inner then we may take $U\in\M(D)$ and then $B=D$.\vs

To obtain a more interesting algebra $B$ we have to assume that $U$ is in a sense independent of elements of $D$. One of the symptom of this independence is the existence of dual action. We say that $B$ admits a dual action if for any $z\in S^{1}$ there exists an automorphism $\beta_{z}$ of $B$ such that $\beta_{z}(U)=zU$ and $\beta_{z}(d)=d$ for any $d\in D$. It turns out that in any case one can find a covariant representation such that the algebra \rf{1744} admits a dual action. Moreover the algebra $B$ with the dual action is unique (up to isomorphism): It does not depend on the choice of covariant representation. This unique  $\cst$-algebra is denoted by $D\ltimes_{\alpha}\inte$ and called the crossed product of $D$ by the automorphism $\alpha$. It is equipped with the dual action $\beta$ of $S^{1}$ and distinguished element $U\in\M(D\ltimes_{\alpha}\inte)$ such that 
\[
\etyk{1700}
\beta_{z}(U)=zU
\]
for any $z\in S^{1}$. \vs

Conversely assume that $B$ is a $\cst$-algebra equipped with an action $\beta$ of $S^{1}$ and a distinguished element $U\in\M(B)$ such that formula \rf{1700} holds. Then
 \[
D=\set{d\in B}{\beta_{z}(d)=d\text{ for all }z\in S^{1}}
\]
is a $\cst$-subalgebra of $B$, $DB=B$, formula \rf {imp} defines an automorphism $\alpha$ of $D$ and relation \rf {1744} holds.\vs

In the above the single automorphism may be replaced by a locally compact quantum group of automorphisms. We shall consider locally compact quantum group $G=(A,\Delta)$.  %In the similar way one defines an left action of $G$ on a $\cst$-algebra.
In this case the formula \rf{imp} takes the form 
\[
\etyk{imp1}
\alpha(d)=U(d\tens\I{A})U^{*},
\]
where $U\in\M(\B_{0}(K)\tens A)$ is a unitary representation of $G$ acting on $K$. It is assumed that $\alpha(d)\in\M(D\tens A)$ for any $d\in D$ and that (using the notation \rf{notacja})
\[
\etyk{0716}
\alpha(D)(\I{D}\tens A)=D\tens A.
\]
Then $\alpha\in\Mor(D,D\tens A)$. We say that $\alpha$ is an action of $G$ on $D$ implemented by a representation $U$. \vs

In general, a right action of $G$ on a $\cst$-algebra $D$ is an injective morphism $\alpha\in\Mor(D,D\tens A)$ such that $(\alpha\tens\id{A})\comp\alpha=(\id{D}\tens\,\Delta)\comp\alpha$. The action is said to be continous if the Podleś condition \rf{0716} holds. In the similar way one defines left actions.
\vs

Any continuous action $\alpha$ of a locally compact group $G$ on a $\cst$-algebra $D$ can be implemented. For any $\alpha$  there exist a Hilbert space $K$ and a pair $(j,U)$, where $j$ is an embedding of $D$ in $\B(K)$ and $U$ is a unitary representation of $G$ acting on $K$ such that (identifying $d$ with $j(d)$) we have \rf{imp1}. We say that $(j,U)$ is a covariant representation of $(D,\alpha)$. Moreover one may assume that $U$ is weakly contained in the regular representation. It means that $U$ is of the form
\[
U=(\psi\tens\id{A})V,
\]
where $\psi\in\Rep(\Ahat,K)$ and $V\in\M(\Ahat\tens A)$ is the canonical bicharacter establishing duality between $G=(A,\Delta)$ and $\Ghat=(\Ahat,\Deltahat)$. In what follows we shall use shorthand\footnote{This is extended  leg numbering notation} notation $V_{\psi2}=(\psi\tens\id{A})V$.
\vs

%Unitary representations of locally compact groups weakly contained in the regular representation are in one-to-one correspondence with representations of their $\cst_{r}$-algebras. Denoting by $\psi\in\Rep(\cst_{r}(G),K)$ the representation corresponding to $U$ we have
%\[
%U_{g}=\psi(V_{g})
%\]
%for any $g\in G$. In this formula $V$ is the canonical embedding of $G$ into the group of unitaries of $\M(\cst_{r}(G))$.

\vs

Again it is interesting to extend $D$ by including $\psi$-copy of $\Ahat$. We shall use notation \rf{notacja} (see next section). Let
\[
\etyk{541}
\begin{array}{r@{\;=\;}l}
B&
\set{\psi(\ahat)d\Vs{3.5}}{\ahat\in\Ahat,\ d\in D}^{\rm CLS}\\
\Vs{5}&\psi(\Ahat)D.
\end{array}
\]
Then $B$ is a $\cst$-algebra, $D,\psi(\Ahat)\subset\M(B)$ and $DB=\psi(\Ahat)B=B$. The latter means that the embedding $D\subset\M(B)$ is a morphism from $D$ into $B$ and $\psi$ is a morphism from $\Ahat$ into $B$. We refer to the next section for the concept of morphism in the category of $\cst$-algebras. In general, for given $D$ and $\alpha$, the algebra $B$ may depend on the used covariant representation. %For instance, if $\alpha$ is inner then we may take $U_{g}\in\M(D)$ for all $g\in G$ and then $B=D$.
\vs

To obtain a more interesting algebra $B$ we have to assume that elements of $\Ahat$ are in a sense independent of elements of $D$. One of the symptom of this independence is the existence of dual action. We say that $B$ admits a dual action if there exists an injective morphism $\beta\in\Mor(B,\Ahat\tens B)$ such that $\beta(\psi(\ahat))=(\id{\Ahat}\tens\,\psi)\Deltahat(\ahat)$ and $\beta(d)=\I{\Ahat}\tens\,d$ for any $d\in D$ and $\ahat\in\Ahat$. It turns out that in any case one can find a covariant representation such that the algebra \rf{541} admits a dual action. Moreover the algebra $B$ with the dual action is unique (up to isomorphism): It does not depend on the choice of covariant representation. This unique $\cst$-algebra is denoted by $D\ltimes_{\alpha}G$ and called the crossed product of $D$ by the action $\alpha$. It is equipped with the dual action $\beta$ of $\Ghat$ and the morphism $\psi\in\Mor(\Ahat,B)$ such that 
\[
\etyk{0747}
\beta(\psi(\ahat))=(\id{\Ahat}\tens\,\psi)\Deltahat(\ahat)
\]
for any $\ahat\in\Ahat$. In general, a triple $(B,\beta,\psi)$, where $B$ is a $\cst$-algebra, $\beta\in\Mor(B,\Ahat\tens B)$ is a left action of $\Ghat$ on $B$ and $\psi$ is an injective morphism from $\Ahat$ into $B$ is called $G$-product if formula \rf{0747} holds.\vs

One of the aims of the Landstad theory was to describe position of $D$ within $\M(D\ltimes_{\alpha}G)$. It could be easily shown that elements $d\in D$ satisfy the following three conditions:\vs

{\rm L1}. $\beta(d)=\I{\Ahat}\tens\,d$,\vs

{\rm L2}. $\psi(\ahat)d\in B$ for any $\ahat\in\Ahat$,\vs

{\rm L3}.  $\left(V_{\psi2}(d\tens\I{A})V_{\psi2}^{*}\Vs{3}\right)(\I{B}\tens a)\in\M(B)\tens A$ for any $a\in A$.\vs

\noindent The above conditions were (in a more classical language) formulated by Magnus Landstad\footnote{see formulae (3.6) - (3.8) in \cite{Lands}. In fact in the second condition Landstad demanded additionally that $d\psi(\ahat)\in B$. However this requirement is redundant, it follows from other Landstad conditions (see last section)} (see also section 7.8 of \cite{Ped}). Assuming that $G$ is a locally compact group Landstad was able to show that his conditions completely characterise those elements of $\M(B)$ that belong to $D$. More than that Landstad proved (for classical $G$) that any $G$-product $(B,\beta,\psi)$ comes from crossed product construction: 
\[
B=D\ltimes_{\alpha}G,
\]
where $D$ is the $\cst$-algebra consisting of all elements of $d\in\M(B)$ satisfying Conditions L1, L2 and L3, and 
\[
\etyk{0835}
\alpha(d)=V_{\psi2}(d\tens\I{A})V_{\psi2}^{*}
\]
for any $d\in D$.
\vs

Instead of looking for conditions characterising elements of $D$ within $\M(D\ltimes_{\alpha}G)$ one may formulate properties concerning the $\cst$-algebra $D$ itself. It could be easily shown that\vs

{\rm V1}. $\beta(d)=\I{\Ahat}\tens\,d$ for any $d\in D$,\vs

{\rm V2}. $B=\psi(\Ahat) D$,\vs

{\rm V3}. $\left(V_{\psi2}(D\tens\I{A})V_{\psi2}^{*}\Vs{3}\right)(\I{B}\tens A)=D\tens A$.\vs

\noindent These conditions were formulated by Stefaan Vaes in \cite{Vaes1}. Assuming that $G$ is a locally compact regular quantum group Vaes was able to show that for any $G$-product $(B,\beta,\psi)$ there exists unique $\cst$-subalgebra $D$ of $\M(B)$ satisfying Conditions V1, V2 and V3. This subalgebra is equipped with the left action $\alpha\in\Mor(D,A\tens D)$ of $G$ introduced by \rf{0835} and $(B,\beta,\psi)$ comes from crossed product construction: $B=D\ltimes_{\alpha}G$. This way Vaes extended Landstad theory to regular locally compact quantum groups.
\vs

It is interesting to compare Landstad and Vaes conditions.  Let $(B,\beta,\psi)$ be a $G$-product, $D\subset\M(B)$ be a $\cst$-subalgebra and $d\in D$. Then L1 coincides with V1, L2 follows from V2 and assuming V3 we see that for any $a\in A$ we have  $\left(V_{\psi2}(d\tens\I{A})V_{\psi2}^{*}\Vs{3}\right)(\I{B}\tens a)\in D\tens A\subset\M(B)\tens A$. It shows that L3 follows from V3. %To understand the importance of Condition V3 one should combine it with \rf{0835}: 

\vs

Let $(B,\beta,\psi)$ be a $G$-product. In general (for non-regular $G$) a subalgebra $D\subset\M(B)$ satisfying Vaes conditions may not exist. To regain the existence of $D$ we have to replace V3 by a weaker Condition C3. Roughly speaking, C3 means that the slices of $V_{\psi2}(D\tens\I{A})V_{\psi2}^{*}$ generate $\cst$-algebra $D$.\vs

In what follows the subalgebra $G\subset\M(B)$ satisfying Conditions V1, V2 and C3 will be called Landstad algebra of $(B,\beta,\psi)$. We shall prove that any $G$-product admits unique Landstad algebra. 
Unfortunately now (when Condition V3 is not satisfied) formula \rf{0835} does not define a action of $G$ on $D$. In general $\alpha(d)\notin\M(D\tens A)$. To deal with the problem we invent the notion of weak action adapted to this situation. In brief instead of \rf{0716} we assume that slices of $\alpha(d)$ belongs to $D$ and that the set of all slices generate $\cst$-algebra $D$. In the following a pair $(D,\alpha)$, where $D$ is a $\cst$-algebra and $\alpha$ is a weak action of $G$ on $D$, will be called a weak $G$-dynamical system. For regular groups the concepts of weak and continuous actions coincides.\vs

Working with weak actions we have to reconsider the concept of crossed product. Given a $G$-dynamical system $(D.\alpha)$, we shall construct a $G$-product $(B,\beta,\psi)$ such that  $D$ plays the role of Landstad algebra of $(B,\beta,\psi)$, $\alpha$ is implemented by $V_{\psi2}$ and $\psi\in\Mor(\Ahat,B)$ is the canonical embedding. In other words $B$ is in a sense crossed product of $D$ by the action $\alpha$ and $\beta$ is the dual action. We are still able to show that $(B,\beta,\psi)$ is unique. It means that the correspondence between $G$-products and $G$-dynamical systems is one to one.\vs

Let us shortly discuss the content of the paper. In section \ref{sek1} we explain the notation used in the paper. In particular we recall the category of $\cst$-algebras (concepts of morphisms and composition of morphisms). We also collect all informations concerning quantum groups used in the paper. Section~\ref{sek2} contains main definitions and results. We introduce (recall) the concepts of $G$-dynamical system and $G$-product and describe the duality between them. In Section \ref{sek3} we analyse the concept of weak action. The most important result is Proposition \ref{HPresults} that enables crossed product construction. At the end we show that for regular groups any weak action is continuous (satisfies Podleś condition). The existence and uniqueness of Landstad algebra for any $G$-product is discussed in section \ref{sek5}. At the end of the section we show the uniqueness of $G$-product corresponding to any weak $G$-dynamical system. Section \ref{sek4} is devoted to the crossed product construction. Next, in  Section \ref{sek8} we discuss in detail weak actions  implemented by a unitary representation of $G$. An application to the Kasprzak version of Rieffel deformation is recalled in Section \ref{sek6}. Finally in Section \ref{sek7} we show (for coameanable $G$) that one of the original Landstad condition is a consequence of the others.

%\newpage

\section{Notation}
\label{sek1}

Throughout the paper we shall use the following notation: For any separable Hilbert space $K$ we set
\[
\begin{array}{r@{\;=\;}l}
\Vs{6}\B(K)&\text{the von Neumann algebra of all bounded operators acting on }K,\\
\Vs{6}\B_{0}(K)&\text{the }\cst\text{-algebra of all compact operators acting on }K,\\
\Vs{6}\B(K)_{*}&\text{the set of all normal functionals on }\B(K)\\
\Vs{5}&\text{the set of all continuous functionals on }\B_{0}(K),\\
\Vs{8}\cst(K)&\set{A\subset\B(K)}{\begin{array}{c}
A\text{ is separable }\cst\text{-algebra}\\ \text{ such that }AK=K
\end{array}}.
\end{array}
\]
\Vs{6}Then $\B(K)$ and $\B_{0}(K)$ are $\cst$-algebras, $\cst(K)$ is a set of $\cst$-algebras, $\B_{0}(K)\in\cst(K)$. In this paper phrase `$\cst$-algebra generated by a set' means `the smallest $\cst$-algebra containing the set'.\vs

We recall that $\B(K)_{*}$ is a bimodule over $\B(K)$. For any $\mu\in\B(K)_{*}$ and $a\in\B(K)$, $\mu a$ and $a\mu$ are normal functionals on $\B(K)$ such that
\[
\etyk{1410}
\begin{array}{r@{\;=\;}l}
(\mu a)(m)&\mu(am),\\
\Vs{5}(a\mu)(m)&\mu(ma)
\end{array}
\]
for any $m\in\B(K)$.\vs

Let $X$ and $Y$ be a norm closed subsets of a $\cst$-algebra. We set
\[
\etyk{notacja}
XY=\set{xy}{
\begin{array}{c}
x\in X\\y\in Y
\end{array}}^{\rm CLS},
\]
where ${\rm CLS}$ stays for norm closed linear span. \vs

For any $\cst$-algebra $A$, $\M(A)$ will denote the multiplier algebra (cf \cite{Ped}) of $A$. Then $A$ is an essential ideal in $\M(A)$. We shall use the category of $\cst$-algebras introduced in \cite{SLW81, Vallin}. It will be denoted by $\cst$. Objects are $\cst$-algebras. For any $\cst$-algebras $A$ and $B$, $\Mor(A,B)$ is the set of all $^{*}$-algebra homomorphisms $\varphi$ from $A$ into $\M(B)$  such that $\varphi(A)B=B$. Any $\varphi\in\Mor(A,B)$ admits unique extension to a unital $^*$-homomorphism $\varphitilde:\M(A)\longrightarrow\M(B)$. If $\varphi\in\Mor(A,B)$ and $\psi\in\Mor(B,C)$ ($A,B,C$ are $\cst$-algebras) then the composition of morphisms $\psi\comp\varphi\in\Mor(A,C)$ is defined as composition of $^{*}$-algebra homomorphisms: $\psi\comp\varphi=\psi\comp\varphitilde$.\vs

Let $A$ be a $\cst$-algebra. Depending on the context we shall use two symbols: $\id{A}$ and $\I{A}$ to denote the identity map acting on $A$. $\id{A}$ will denote the identity morphism acting on $A$, whereas $\I{A}$ will denote the unit element of the multiplier algebra $\M(A)$\footnote{Identifying elements of $\M(A)$ with left multipliers acting on $A$ we have $\I{A}=\id{A}$.}. To simplify notation we write $\I{K}$ and $\id{K}$ instead of $\I{B_{0}(K)}$ and $\id{B_{0}(K)}$. If $A\in\cst(K)$ then $\I{A}=\I{K}$. We shall often omit the index `$A$'  (or `$K$') when the algebra $A$ (or the Hilbert space $K$) is obviously implied by the context.\vs

By definition, representations acting on a Hilbert space $K$ are morphisms into $\B_{0}(K)$: for any $\cst$-algebra $X$, $\Rep(X,K)=\Mor(X,\B_{0}(K))$. We know that $\M(\B_{0}(K))=\B(K)$. Therefore representations are non-degenerate $^{*}$-algebra homomorphisms into $B(K)$. $^{*}$-algebra homomorphism $\pi:X\longrightarrow\B(K)$ is non-degenerate if $\pi(X)\B_{0}(K)=\B_{0}(K)$. Equivalently $\pi$ is non-degenerate if $0\in K$ is the only vector killed by $\pi(x)$ for all $x\in X$.\vs

The category $\cst$ is equipped with a monoidal structure. For any $\cst$-algebras $X$ and $Y$, the tensor product $X\tens Y$ is a $\cst$-algebra. One can easily define tensor product of morphisms. Then $\tens$ becomes an associative functor from $\cst\times\cst$ to $\cst$. In this paper we use exclusively minimal (spacial) tensor product of $\cst$-algebras.\vs

\begin{Prop}
\label{XYYX}
Let $X,Y\in\cst(K)$, where $K$ is a Hilbert space. Assume that
\[
XY=YX.
\]
Then $Z=XY\in\cst(K)$, $X,Y\subset\M(Z)$ and the embeddings are morphisms from $X$ and $Y$ into $Z$. 
\end{Prop}
We shall refer to this situation by saying that $Z$ is a crossed product of $X$ and $Y$. The obvious proof is left to the reader.\vs

We shall use the following shorthand notation borrowed from \cite{BSV}: If $\left(X_{\omega}\right)_{\omega\in\Omega}$ is a family of subsets of a Banach space $X$, then the smallest Banach subspace of $X$ containing all $X_{\omega}$ ($\omega\in\Omega$) will be denoted by
\[
\Vset{X_{\omega}}{\Vs{4}\omega\in\Omega}.
\]
So we have:
\[
\Vset{X_{\omega}}{\Vs{4}\omega\in\Omega}=\left(\;\bigcup_{\omega\in\Omega}X_{\omega}\right)^{\rm CLS}.
\]
Typically we shall deal with expressions of the form 
\[
\Vset{(\omega\tens\id{K})Z}{\Vs{4}\omega\in\B(H)_{*}},
\]
where $Z$ is a linear subset of $\B(H\tens K)$ and $H$ and $K$ are Hilbert spaces. Let $A\in\cst(H)$. By the factorisation theorem \cite{factorization} any $\omega\in\B(H)_{*}$ is of the form $\omega'a$ and of the form $a\omega'$ where $\omega'\in\B(H)_{*}$ and $a\in A$. Therefore
\[
\etyk{facto}
\begin{array}{r@{\;=\;}l}
\Vset{(\omega\tens\id{K})Z}{\Vs{4}\omega\in\B(H)_{*}}
&\Vset{(\omega\tens\id{K})(A\tens\I{K})Z}{\Vs{4}\omega\in\B(H)_{*}}\\
\Vs{7}&\Vset{(\omega\tens\id{K})Z(A\tens\I{K})}{\Vs{4}\omega\in\B(H)_{*}}.
\end{array}
\]
We shall refer to these formulas by saying that $\omega$ {\it emits} $A$ to the right (upper formula) or to the left (lower formula). Reading from the right%hand side 
\ we say that $\omega$ {\it absorbs} $A$.\vs

In the present paper we consider actions of a locally compact quantum group on $\cst$-algebras. The group will be denoted by $G$. We shall not use the full power of the theory developed by Kustermans and Vaes in \cite{Vaes} (see also \cite{SLW03}). Instead we shall assume that $G$ is constructed from a manageable (modular) multiplicative unitary $V$ in the way described in \cite{SLW96c} (see also \cite{SLW01c, SLW07}). In particular we do not assume the existence of Haar measures. \vs
 
 The Haar measure plays an essential role in the Landstad considerations. To construct $\beta$ invariant elements in $\M(B)$ he applies $\beta_{g}$ ($g\in G$) to specially chosen elements of $B$ and then integrates over $G$ using the Haar measure. On the other hand in Vaes approach \cite{Vaes1} (at least in the part that intersects with our interest) the Haar measure plays purely decorative role and can be removed from considerations. Then the proofs become more transparent. \vs
 
 In this paper a locally compact quantum group $G$ is a pair $G=(A,\Delta)$, where $A$ is a $\cst$-algebra and $\Delta\in\Mor(A,A\tens A)$ having a number of properties. One of the property is coassociativity of $\Delta$:
 \[
 (\Delta\tens\id{A})\comp\Delta= (\id{A}\tens\,\Delta)\comp\Delta,
 \]
another one is the cancelation property:
\[
(A\tens\I{A})\Delta(A)=A\tens A=\Delta(A)(\I{A}\tens A).
\]
 
Instead of listing the properties we assume that $(A,\Delta)$ comes from a manageable multiplicative unitary by the construction described in \cite{SLW96c}. Elements of $A$ should be considered as {\it continuous vanishing at infinity functions} on $G$, whereas $\Delta$ encodes the {\it group multiplication} on $G$. Locally compact quantum groups appear in pairs. For any $G=(A.\Delta)$ we have dual group $\Ghat=(\Ahat,\Deltahat)$. The duality between $G$ and $\Ghat$ is described by a bicharacter $V\in\M(\Ahat\tens A)$. It satisfies bicharacter equations:
\[
\etyk{1037}
\begin{array}{r@{\;=\;}l}
(\id{\Ahat}\tens\,\Delta)V&V_{12}V_{13},\\
\Vs{5}(\Deltahat\tens\id{A})V&V_{23}V_{13},
\end{array}
\]
The role of $G$ and $\Ghat$ is symmetric. Replacing $V$ by $\Vhat=\flip(V^{*})$ we obtain bicharacter describing duality between $\Ghat$ and $G$.\vs

An important role in the theory of quantum groups is played by Heisenberg and anti-Heisenberg pairs. These are pairs of representations of $A$ and $\Ahat$ acting on the same Hillbert space. Let $H$ be a Hilbert space, $\sigma\in\Rep(A,H)$ and $\sigmahat\in\Rep(\Ahat,H)$. We say that $(\sigma,\sigmahat)$ is a Heisenberg pair if
\[
\etyk{Ha}
V_{\sigmahat3}V_{1\sigma}=V_{1\sigma}V_{13}V_{\sigmahat3}.
\]
Similarly let $\Hbar$ be a Hilbert space, $\rho\in\Rep(A,\Hbar)$ and $\rhohat\in\Rep(\Ahat,\Hbar)$. We say that $(\rho,\rhohat)$ is an anti-Heisenberg pair if
\[
\etyk{aHa}
V_{1\rho}V_{\rhohat3}=V_{\rhohat3}V_{13}V_{1\rho}.
\]
It is known that representations appearing in Heisenberg and in anti-Heisenberg pairs are faithful.\vs

In the above formulas extended leg numbering notation is used. Both sides of \rf{Ha} belongs to $\M(\Ahat\tens\B_{0}(H)\tens A)$. By definition $V_{\sigmahat3}=\I{\Ahat}\tens(\sigmahat\tens\id{A})V$ and $V_{1\sigma}=(\id{\Ahat}\tens\sigma)V\tens\I{A}$. Similarly both sides of \rf{aHa} belongs to $\M(\Ahat\tens\B_{0}(\Hbar)\tens A)$,  $V_{1\rho}=(\id{\Ahat}\tens\rho)V\tens\I{A}$ and $V_{\rhohat3}=\I{\Ahat}\tens(\rhohat\tens\id{A})V$.\vs

 Inserting on \rf{Ha} $V=\flip(\Vhat^{*})$ one can easily show that $\Vhat_{\sigma3}\Vhat_{1\sigmahat}=\Vhat_{1\sigmahat}\Vhat_{13}\Vhat_{\sigma3}$. It shows that $(\sigmahat,\sigma)$ is a Heisenberg pair for $\Ghat$.\vs

The existence of Heisenberg and anti-Heisenberg pairs is one of the basic features of the theory of quantum groups. In \cite{SLW96c} we start with a multiplicative unitary operator acting on $H\tens H$. Then $A$ and $\Ahat$ appear as $\cst$-algebras acting on $H$ and the embeddings $A\hookrightarrow\B(H)$ and $\Ahat\hookrightarrow\B(H)$ form a Heisenberg pair. Let $(\sigma,\sigmahat)$ be a Heisenberg pair acting on a Hilbert space $H$. Using unitary antipodes $R$ and $\Rhat$  (see \cite{SLW96c}) we can construct anti-Heisenberg pair acting on $\Hbar$. $\Hbar$ is the Hilbert space complex-conjugate to $H$. Then we have canonical anti-unitary mapping $H\ni x\mapsto\xbar\in\Hbar$ and transposition map $\B(H)\ni m\mapsto m^{\top}\in\B(\Hbar)$ introduced by the formula $m^{\top}\xbar=\overline{m^{*}x}$ for any $x\in H$. For any $a\in A$ and $\ahat\in\Ahat$ we set
\[
\begin{array}{r@{\;=\;}l}
\rho(a)&\sigma(a^{R})^{\top},\\
\rhohat(\ahat)\Vs{4}&\sigmahat(\ahat^{\Rhat})^{\top}.
\end{array}
\]
Then $(\rho,\rhohat)$ is an anti-Heisenberg pair acting on $\Hbar$.\vs

Let $(\sigma,\sigmahat)$ be a Heisenberg pair acting on $H$ and $M$ ($\Mhat$ resp.) be the weak closure of $\sigma(A)$ ($\sigmahat(\Ahat)$ resp.). We denote by $A_{*}$ ($\Ahat_{*}$ resp.) the set of all linear functionals on $A$ ($\Ahat$ resp.) that admit extensions to normal functionals on $M$ ($\Mhat$ resp):
\[
\begin{array}{r@{\;=\;}l}
A_{*}&\set{\Vs{3}\mu\comp\sigma}{\mu\in\B(H)_{*}},\\
\Ahat_{*}\Vs{5}&\set{\Vs{3}\mu\comp\sigmahat}{\mu\in\B(H)_{*}}.
\end{array}
\]
We know that $\sigma$ and $\sigmahat$ are faithful. Therefore $A_{*}$ and $\Ahat_{*}$ are weakly dense in the set of all continuous functionals on $A$ and $\Ahat$ resp.. It turns out that $A_{*}$ and $\Ahat_{*}$ are independent of the choice of Heisenberg pair. We may even replace $(\sigma,\sigmahat)$ by an anti-Heisenberg pair.\vs

Elements of $A_{*}$ are called normal functionals on $A$. A representation $\phi\in\Rep(A,K)$ is called normal if its matrix elements are normal: $\mu\comp\phi\in A_{*}$ for any $\mu\in\B(K)_{*}$. Similarly one defines normal functionals on $\Ahat$ and normal representations of $\Ahat$. Representations in Heisenberg and anti-Heisenberg pairs are normal.\vs

Algebras $A$ and $\Ahat$ coincide with the closures of the sets of slices of $V$:
\[
\etyk{algabry}
\begin{array}{r@{\;=\;}l}
A&\set{(\omegahat\tens\id{A})V}{\omegahat\in\Ahat_{*}}^{\rm norm\ closure},\\
\Ahat\Vs{6}&\set{(\id{\Ahat}\tens\,\omega)V}{\Vs{4}\omega\in A_{*}}^{\rm norm\ closure}.
\end{array}
\]

Let $(\sigma,\sigmahat)$ be a Heisenberg pair and $(\rho,\rhohat)$ be an anti-Heisenberg pair. Combining \rf{1037} with \rf{Ha} and \rf{aHa} we get
\[
\begin{array}{r@{\;=\;}l}
(\id{\Ahat}\tens\,(\sigma\tens\id{A})\Delta)V&V_{1\sigma}V_{13}=V_{\sigmahat3}V_{1\sigma}V_{\sigmahat3}^{*},\\
\Vs{5}((\id{\Ahat}\tens\rhohat)\Deltahat\tens\id{A})V&V_{\rhohat3}V_{13}=V_{1\rho}V_{\rhohat3}V_{1\rho}^{*},\\
\Vs{5}\left(\id{\Ahat}\tens\,(\rho\tens\id{A})\comp\flip\comp\Delta\right)V&
V_{13}V_{1\rho}=V\Vs{3}_{\rhohat3}^{*}V_{1\rho}V_{\rhohat3},\\
\Vs{5}\left((\id{\Ahat}\tens\,\sigmahat)\comp\flip\comp\Deltahat\tens\id{A}\right)V&
V_{13}V_{\sigmahat3}=
V\Vs{3}_{1\sigma}^{*}V_{\sigmahat3}V_{1\sigma}
\end{array}
\]
Computing appropriate slices and taking into account \rf{algabry} we obtain
\[
\etyk{impleme}
\begin{array}{r@{\;=\;}l}
(\sigma\tens\id{A})\Delta(a)&V_{\sigmahat2}(\sigma(a)\tens\I{A})V_{\sigmahat2}^{*},\\
(\id{\Ahat}\tens\rhohat)\Deltahat(\ahat)\Vs{5}&V_{1\rho}(\I{\Ahat}\tens\rhohat(\ahat))V_{1\rho}^{*},\\
(\rho\tens\id{A})\comp\flip\comp\Delta(a)\Vs{5}
&V\Vs{3}_{\rhohat2}^{*}\left(\rho(a)\tens\I{A}\right)V_{\rhohat2},\\
(\id{\Ahat}\tens\,\sigmahat)\comp\flip\comp\Deltahat(\ahat)\Vs{5}
&V\Vs{3}_{1\sigma}^{*}\left(\I{\Ahat}\tens\,\sigmahat(\ahat)\right)V_{1\sigma}
\end{array}
\]
for any $a\in A$ and $\ahat\in\Ahat$.\vs

We shall use convolution product of normal functionals. For any $\omega',\omega\in A_{*}$ we set 
 \[
 \etyk{2241}
 \omega'*\omega=(\omega'\tens\,\omega)\comp\Delta.
 \]
Using first formula of \rf{impleme} one can easily show that $\omega'*\omega\in A_{*}$. \vs

%\SLW{
%Let $(\sigma,\sigmahat)$ be a Heisenberg pair acting on a Hilbert space $H$. We say that $(\sigma,\sigmahat)$ is canonical if there exists a unitary operator $T\in\M(A\tens\B_{0}(H))$ such that
%\[
%\begin{array}{r@{\;=\;}l}
%(\id{A}\tens\,\sigma)\Delta(a)&T^{*}(\I{A}\tens\,\sigma(a))T\\
%\Vs{4}\I{A}\tens\,\sigmahat(\ahat)&T^{*}(\I{A}\tens\,\sigmahat(\ahat))T
%\end{array}
%\]
%for any $a\in A$ and $\ahat\in\Ahat$. Equivalently, $(\sigma,\sigmahat)$ is canonical if there exists a unitary $U\in\M(\Ahat\tens\B_{0}(H))$ such that
%\begin{align}
%\I{\Ahat}\tens\,\sigma(a)&=U(\I{\Ahat}\tens\,\sigma(a))U^{*}\label{F1532}\\
%\Vs{4}(\id{\Ahat}\tens\,\sigmahat)\Deltahat(\ahat)&=U(\I{\Ahat}\tens\,\sigmahat(\ahat))U^{*}\label{S1532}
%\end{align}
%for any $a\in A$ and $\ahat\in\Ahat$. We shall use a canonical Heisenberg pair in section \ref{sek4} to construct crossed product. To remove possible doubts about this construction we have to show that there exists a canonical Heisenberg pair.  }

Let $(\sigma,\sigmahat)$ be a Heisenberg pair acting on a Hilbert space $H$. We say that $(\sigma,\sigmahat)$ is canonical if there exists $\rho\in\Rep(A,H)$ such that $(\rho,\sigma)$ is a commuting pair\footnote{it means that $[\rho(a),\sigma(a')]=0$ for any $a,a'\in A$} and $(\rho,\sigmahat)$ is an anti-Heisenberg pair. We shall use a canonical Heisenberg pair in section \ref{sek4} to construct crossed product. To remove possible doubts about this construction we have to show that there exists a canonical Heisenberg pair.

\begin{Prop}
\label{canonHeis}
Let $(\sigma,\sigmahat)$ be a Heisenberg pair acting on a Hilbert space $H$ and $(\rho,\rhohat)$ be an anti- Heisenberg pair acting on a Hilbert space $\Hbar$ . For any $a\in A$ and $\ahat\in\Ahat$ we set
\[
\begin{array}{r@{\;=\;}l}
\sigma'(a)&\I{\Hbar}\tens\,\sigma(a),\\
\Vs{5}\sigmahat'(\ahat)&(\rhohat\tens\,\sigmahat)\Deltahat(\ahat).
\end{array}
\]
Then $(\sigma',\sigmahat')$ is a canonical Heisenberg pair acting on $H'=\Hbar\tens H$.
\end{Prop}
\begin{pf} 
We shall prove that $(\sigma',\sigmahat')$ is a Heisenberg pair. 
In the following considerations we deal with elements of $\M(\Ahat\tens\B_{0}(\Hbar\tens H)\tens A)=\M(\Ahat\tens\B_{0}(\Hbar)\tens\B_{0}(H)\tens A)$. Clearly $V_{1\sigma'}=V_{1\sigma}$ and $V_{\sigmahat'4}=V_{\sigmahat4}V_{\rhohat4}$ (the latter follows from the second formula of \rf{1037}). We assumed that $(\sigma,\sigmahat)$ is a Heisenberg pair: $V_{\sigmahat4}V_{1\sigma}=V_{1\sigma}V_{14}V_{\sigmahat4}$. The reader should notice that {\it nontrivial legs} of $V_{1\sigma}$ and $V_{\rhohat4}$ belong to disjoint sets. Therefore the two unitaries commute. Now we have:
\[
V_{\sigmahat'4}V_{1\sigma'}
=V_{\sigmahat4}V_{\rhohat4}V_{1\sigma}
=V_{\sigmahat4}V_{1\sigma}V_{\rhohat4}
=V_{1\sigma}V_{14}V_{\sigmahat4}V_{\rhohat4}
=V_{1\sigma'}V_{14}V_{\sigmahat'4}
\]
It shows that $(\sigma',\sigmahat')$ is a Heisenberg pair.\vs

For any $a\in A$ we set $\rho'(a)=\rho(a)\tens\I{H}$. Then $\rho'\in\Rep(A,\Hbar\tens H)$. 
We shall prove that $(\rho',\sigmahat')$ is an anti-Heisenberg pair.  Clearly $V_{1\rho'}=V_{1\rho}$ and $V_{\sigmahat'4}=V_{\sigmahat4}V_{\rhohat4}$. We assumed that $(\rho,\rhohat)$ is an anti-Heisenberg pair: $V_{1\rho}V_{\rhohat4}=V_{\rhohat4}V_{14}V_{1\rho}$. The reader should notice that {\it nontrivial legs} of $V_{1\rho}$ and $V_{\sigmahat4}$ belong to disjoint sets. Therefore the two unitaries commute. Now we have:
\[
V_{1\rho'}V_{\sigmahat'4}
=V_{1\rho}V_{\sigmahat4}V_{\rhohat4}
=V_{\sigmahat4}V_{1\rho}V_{\rhohat4}
=V_{\sigmahat4}V_{\rhohat4}V_{14}V_{1\rho}
=V_{\sigmahat'4}V_{14}V_{1\rho'}.
\]
It shows that $(\rho',\sigmahat')$ is an anti-Heisenberg pair. Obviously $(\rho'.\sigma')$ is a commuting pair. Therefore the Heisenberg pair $(\sigma',\sigmahat')$ is canonical.
\end{pf}\vs

An important concept of the theory of locally compact quantum groups is that of regularity and semi-regularity. It was introduced by Baaj and Skandalis in \cite{bs} and Baaj in \cite{BaajE2}. Let $(\sigma,\sigmahat)$ be a Heisenberg pair acting on $H$, $W=V_{\sigmahat\sigma}=(\sigmahat\tens\,\sigma)V$ and%$V$ be a manageable multiplicative unitary and
\[
C=\set{(\id{H}\tens\,\omega)(\Sigma W)}{\Vs{4}\omega\in\B(H)_{*}}^{\rm CLS},
\]
where $\Sigma$ be the flip operator acting on $H\tens H$: $\Sigma(x\tens y)=y\tens x$ for all $x,y\in H$.
We say that $W$ is semi-regular if $\B_{0}(H)\subset C$ and regular if $\B_{0}(H)=C$. %It turns out that regularity (and semiregularity) is a property of the corresponding locally compact quantum group. Indeed, 
It was shown in \cite{bs} that $W$ is regular if and only if
\[
\etyk{reg}
(\Ahat\tens\I{A})V(\I{\Ahat}\tens A)=\Ahat\tens A.
\]
This formula no longer depends on the choice of Heisenberg pair. This is an equality of subsets of $\M(\Ahat\tens A)$. Due to this fact regularity is the property of the group (not of a particular multiplicative unitary $W$ related to the group). It is known that $V^{\Rhat\tens R}=V$. Applying anti-multiplicative involution $\Rhat\tens R$ to the both sides \rf{reg} we obtain equivalent condition:
\[
\etyk{reg1}
(\I{\Ahat}\tens A)V(\Ahat\tens\I{A})=\Ahat\tens A.
\]
%One can easily show that the dual of a regular group is regular. Therefore \rf{reg1} is equivalent to
%\[
%\etyk{reg2}
%(\I{A}\tens\Ahat)\Vhat(A\tens\I{\Ahat})=A\tens\Ahat.
%\]

%\newpage

\section{Main definitions and results} 
\label{sek2}

Let $G=(A,\Delta)$ be a locally compact quantum group and $D$ be a $\cst$-algebra. A weak action of $G$ on $D$ is a bilinear mapping $\alpha:A_{*}\times D\longrightarrow D$ subject to a number of conditions. One of the conditions is faithfulness: $\alpha$ is called faithful if for any non-zero $d\in D$ there exists $\omega\in A_{*}$ such that $\alpha(\omega,d)\neq 0$.\vs

To formulate the definition of weak action in a compact form we need the concept of faithful presentation.\vs

\begin{Def}
Let $K$ and $H$ be Hilbert spaces, $\pi\in\Rep(D,K)$ and $\pitilde\in\Mor(D,\B_{0}(K)\tens A)$. We say that $(\pi,\pitilde)$ is a faithful presentation of $\alpha$ on $K$ if $\pi$ is faithful and
\[
\etyk{0855}
\pi\left(\alpha(\omega,d)\Vs{3.5}\right)=(\id{K}\tens\,\omega)\pitilde(d)
\]
for any $d\in D$ and $\omega\in A_{*}$. 
\end{Def}
Clearly $\pitilde$ is injective for faithful $\alpha$. Now we are ready to formulate our main definition.

%
%\begin{Def}
%\label{multipli}
%Let $D$ be a $\cst$-algebra and $\alpha$ be a bilinear mapping acting from $A_{*}\times D$ into $D$. We say that $\alpha$ is multiplicative if for any positive $\omega\in A_{*}$ and any sequence $\varphi_{0},\varphi_{1},\varphi_{2},\dots\in A_{*}$ such that for all $a\in A$
%\[
%\omega(a^{*}a)=\sum_{k=0}^{\infty}\modul{\varphi_{k}(a)}^{2},
%\]
%we have 
%\[
%\etyk{2248}
%\alpha(\omega,d^{*}d)=\sum_{k=0}^{\infty}\alpha(\varphi_{k},d)^{*}\alpha(\varphi_{k},d),
%\]
%where $d$ runs over $D$ and the series is norm convergent.
%\end{Def}

\begin{Def}
\label{mea}
Let $D$ be a $\cst$-algebra and $\alpha$ be a faithful bilinear mapping from $A_{*}\times D$ into $D$. We say that $\alpha$ is a weak right action of $G$ on $D$ if the following three conditions are fulfilled:\vs

{\rm WA0.} $\alpha$ admits a faithful presentation.\vs

{\rm WA1.} $D$ is generated by $\Vset{\alpha(\omega,D)}{\Vs{3.5}\omega\in A_{*}}$,\vs

{\rm WA2.} $\alpha(\omega,\alpha(\omega',d))=\alpha(\omega*\omega',d)$ for any $d\in D$ and $\omega,\omega'\in A_{*}$.
\end{Def}\vs

Alternatively we say that the pair $(D,\alpha)$ is a weak $G$-dynamical system. Let $(D,\alpha)$ and $(D',\alpha')$ be weak $G$-dynamical systems and $\imath:D\rightarrow D'$ be an isomorphism of $\cst$-algebras. We say that $\imath$ is an isomorphism of dynamical systems if $\imath$ intertwins $\alpha$ and $\alpha'$:
\[
\imath(\alpha(\omega,d))=\alpha'(\omega,\imath(d))
\]
for any $d\in D$ and $\omega\in A_{*}$.\vs

It would be nice to replace Condition WA0 by a number of algebraic and topological conditions imposed directly on the bilinear map $\alpha$. This issue goes beyond the present paper. Eventually we shall return to it later. \SLW{}\vs

Let $(D,\alpha)$ be a weak $G$-dynamical system and
\[
\etyk{910}
D_{1}=\Vset{\alpha(\omega,D)}{\Vs{3.5}\omega\in A_{*}}.
\]
We say that $D_{1}$ is a range of $\alpha$. The range $D_{1}$ is a norm closed subspace of $\cst$-algebra $D$, so it carries a structure of an operator space in the sense of \cite{EffRuan}. Moreover $D_{1}$ is invariant under Hermitian conjugation. In general (when $G$ is not semiregular) $D_{1}$ is not a $\cst$-algebra (see \cite[Proposition 5.7]{BSV}). According to Condition WA1, $D_{1}$ generates the $\cst$-algebra $D$. Therefore
\[
D=\Vset{D_{k}}{\Vs{4}k=1,2,3,\dots},
\]
where  $D_{2}=D_{1}D_{1}$, $D_{3}=D_{1}D_{1}D_{1}$ and so on. Let $(\pi,\pitilde)$ be a faithful presentation of $\alpha$. Then \rf{0855} shows that
\[
\etyk{1248}
\pi(D_{1})=\Vset{(\id{K}\tens\,\omega)
\pitilde(D)}{\Vs{4}\omega\in A_{*}}.
\]
\vs

The weakness is a weak condition imposed on actions of quantum groups on $\cst$-algebras. Stronger are weak continuity and continuity \cite{BSV}. We say that $\alpha$ is a weakly continuous action of $G$ if there exixts $\alphatilde\in\Mor(D,D\tens A)$ such that
\[
\alpha(\omega,d)=(\id{D}\tens\,\omega)\alphatilde(d)
\]
for any $\omega\in A_{*}$ and $d\in D$. If moreover the Podleś condition $\alphatilde(D)(\I{D}\tens A)=D\tens A$ is satisfied then $\alpha$ is called continuous. Clearly continuous actions are weakly continuous.\vs

\begin{Thm}
\label{BSV}
Let $\alpha$ be a weak action of a regular locally compact quantum group $G$ on a $\cst$-algebra $D$. Then $\alpha$ is a continuous action.
\end{Thm}

This theorem belongs essentially to Baaj, Skandalis and Vaes \cite{BSV}. We shall present a proof because our setting is slightly more general than the one used in \cite{BSV} and because we deal with right actions. Left actions considered in \cite{BSV} did not require the use of the anti-Heisenberg pair.

\begin{Def}
\label{GP}
$G$-product is a triple $(B,\beta,\psi)$, consisting of a $\cst$-algebra $B$, a left continuous action $\beta\in\Mor(B,\Ahat\tens B)$ of $\Ghat$ on $B$ and an injective morphism $\psi\in\Mor(\Ahat,B)$ such that the diagram
\[
\etyk{a1}
\vcenter{
\xymatrix{
\Ahat\ar[rr]^-{\psi}\ar[d]_-{\Deltahat}&&B\ar[d]^-{\beta}\\
\Ahat\tens\Ahat\ar[rr]_-{\id{\Ahat}\tens\,\psi}&&\Ahat\tens B
}}
\]
is commutative.
\end{Def}

We recall that $\beta$ is a continuous left action of $\Ghat$ on $B$ if $\beta$ is an injective morphism $\beta\in\Mor(B,\Ahat\tens B)$ such that $(\id{\Ahat}\tens\beta)\beta=(\Deltahat\tens\id{B})\beta$ and 
\[
\etyk{1733}
(\Ahat\tens\I{B})\beta(B)=\Ahat\tens B.
\]

Let $(B^{1},\beta^{1},\psi^{1})$ and $(B^{2},\beta^{2},\psi^{2})$ be  $G$-products and $\jmath:B^{1}\longrightarrow B^{2}$ be an isomorphism of $\cst$-algebras. We say that $\jmath$ is an isomorphism of the $G$-products if the diagram
\[
\etyk{a2}
\vcenter{
\xymatrix{
&&B^{1}\ar[rr]^-{\beta^{1}}\ar[dd]^-{\jmath}&&\Ahat\tens B^{1}\ar[dd]^-{\id{\Ahat}\tens\jmath}\\
\Ahat\ar[rru]^-{\psi^{1}}\ar[rrd]_-{\psi^{2}}&&&&\\
&&B^{2}\ar[rr]_-{\beta^{2}}&&\Ahat\tens B^{2}
}}
\]
is commutative.

\begin{Def}
\label{LC}
Let  $(B,\beta,\psi)$ be a $G$-product and $D\subset\M(B)$ be a $\cst$-subalgebra. We say that $D$ is a Landstad algebra of $(B,\beta,\psi)$  if the following three  conditions are satisfied:\vs

{\rm C1}. $\beta(d)=\I{\Ahat}\tens\,d$ for any $d\in D$,\vs

{\rm C2}. $B=\psi(\Ahat) D$,\vs

{\rm C3}. The $\cst$-algebra generated by 
\[
\Vset{(\id{B}\tens\,\omega)\left(V_{\psi2}(D\tens\I{H})V_{\psi2}^{*}\right)}{\omega\in A_{*}}
\]
coincides with $D$.
\end{Def}
We know that $B\psi(\Ahat)=B$ (this is because $\psi\in\Mor(\Ahat,B))$. Therefore $BD=B\psi(\Ahat)D=BB=B$. It shows that Landstad algebra is a non-degenerate subalgebra od $\M(B)$. In other words the embedding $D\hookrightarrow\M(B)$ is a morphism from $D$ into $B$ (it belongs to $\Mor(D,B)$). Let $(B,\beta,\psi)$ be a $G$-product and $D$ be its Landstad algebra. Then using Condition C1 and commutativity of \rf{a2} we get
\[
\etyk{2346}
\beta\left(\psi(\ahat)\Vs{3}d\right)=\left((\id{\Ahat}\tens\,\psi)\Deltahat(\ahat)\right)(\I{\Ahat}\tens\,d).
\]
for any $\ahat\in\Ahat$ and $d\in D$. According to Condition C2, this formula determines $\beta$ completely.\vs

\begin{Thm}
\label{Lan}
Any $G$-product admits one and only one Landstad algebra. Let $(B,\beta,\psi)$ be a $G$-product and $D$ be its Landstad algebra. Then the formula
\[
\etyk{1903}
\alpha(\omega,d)=(\id{B}\tens\,\omega)\left(V_{\psi2}(d\tens\I{A})V_{\psi2}^{*}\right)
\] 
for any $\omega\in A_{*}$ and $d\in D$ defines a weak right action $\alpha$ of $G$ on $D$.
\end{Thm}
The resulting weak $G$-dynamical system $(D,\alpha)$ will be called Landstad $G$-dynamical system of $(B,\beta,\psi)$. It turns out that Landstad $G$-dynamical system determines the related $G$-product.
\begin{Thm}
\label{uniqueness}
Let $(B_{1},\beta_{1},\psi_{1})$ and $(B_{2},\beta_{2},\psi_{2})$ be $G$-products, $(D_{1},\alpha_{1})$ and $(D_{2},\alpha_{2})$ be corresponding Landstad $G$-dynamical systems and $\imath:D_{1}\longrightarrow D_{2}$ be an isomorphism of $(D_{1},\alpha_{1})$ and $(D_{2},\alpha_{2})$. Then $G$-products $(B_{1},\beta_{1},\psi_{1})$ and $(B_{2},\beta_{2},\psi_{2})$ are isomorphic: there exists unique isomorphism $\jmath:B_{1}\longrightarrow B_{2}$ of the $G$-products extending $\imath$.
\end{Thm}

To establish the one to one correspondence between weak $G$-dynamical systems and $G$-products we have to construct a $G$-product corresponding to a given weak $G$-dynamical system. The construction is described in the following theorem. 

\begin{Thm}%[Crossed product construction]
\label{CP}

Let $(D,\alpha)$ be a weak $G$-dynamical system, $(\pi,\pitilde)$ be a faithful presentation of $\alpha$ on a Hilbert space $K$ and $(\sigma,\sigmahat)$ be a canonical Heisenberg pair  acting on a Hilbert space $H$.  For any $\ahat\in\Ahat$ we set $\psi(\ahat)=\I{K}\tens\,\sigmahat(\ahat)$. Let $B=\psi(\Ahat)(\id{K}\tens\,\sigma)\pitilde(D)$. Then $B\subset\B(K\tens H)$ and \vs

{\rm 1.} $B$ is a $\cst$-algebra, $\psi(\Ahat),\:(\id{\Ahat}\tens\,\sigma)\pitilde(D)\subset\M(B)$ and the embeddings are morphisms from $\psi(\Ahat)$ and $(\id{K}\tens\,\sigma)\pitilde(D)$ into $B$. Consequently $\psi\in\Mor(\Ahat,B)$ and $(\id{K}\tens\,\sigma)\pitilde\in\Mor(D,B)$.\vs

{\rm 2.} $B$ is equipped with unique continuous left action $\beta\in\Mor(B,\Ahat\tens B)$ of $\Ghat$ such that for any $\ahat\in\Ahat$ and $d\in D$ we have
\[
\begin{array}{r@{\;=\;}l}
\beta\left(\psi(\ahat)\Vs{3}\right)&(\id{\Ahat}\tens\,\psi)\Deltahat(\ahat),\\
\Vs{5}\beta\left(\Vs{3}(\id{K}\tens\,\sigma)\pitilde(d)\right)&\I{\Ahat}\tens\,(\id{K}\tens\,\sigma)\pitilde(d).
\end{array}
\]

{\rm 3.} $(B,\beta,\psi)$ is a $G$-product,  Landstad algebra of $(B,\beta,\psi)$ coincides with $(\id{K}\tens\,\sigma)\pitilde(D)$. Original weak $G$-dynamical system $(D,\alpha)$ is isomorphic to the Landstad $G$-dynamical system of $(B,\beta,\psi)$. More precisely the map $(\id{K}\tens\,\sigma)\pitilde:D\longrightarrow(\id{K}\tens\,\sigma)\pitilde(D)$ is an isomorphism of the weak $G$-dynamical systems.
\end{Thm}\vs

The $\cst$-algebra $B$ appearing in the above Theorem is called a crossed product of $D$ by the action $\alpha$ of $\Ghat$ and denoted by $D\rtimes_{\alpha}\Ghat$. In the same case $\beta$ is called the dual action and often denoted by $\alphahat$.%\newpage

\section{Weak actions}
\label{sek3}

In this section we investigate properties of $G$-dynamical systems. In particular we shall prove Theorem \ref{BSV}. We start with the following Proposition describing analytical properties that follow from Definition \ref{mea}. 
\begin{Prop}
\label{WA2}
Let $(D,\alpha)$ be a $G$-dynamical system, $(\pi,\pitilde)$ be a faithful presentation of $\alpha$ on a Hilbert space $K$, $(\sigma,\sigmahat)$ be a Heisenberg pair acting on a Hilbert space $H$ and $(\rho,\rhohat)$ be an anti-Heisenberg pair acting on a Hilbert space $\Hbar$. Then for any $\omega\in A_{*}$, $\mu\in\B(\Hbar)_{*}$ and $d\in D$ we have
\begin{align}
\pitilde(\alpha(\omega,d))_{1\sigma}=&\ 
(\id{K\tens H}\tens\,\omega)
\left(V_{\sigmahat3}\pitilde(d)_{1\sigma}V\Vs{3}_{\sigmahat3}^{*}\right),\label{1708}
\\
\pitilde(\alpha(\mu\comp\rho,d))=&\ 
(\id{K}\tens\,\mu\tens\id{A})
\left(V\Vs{3}_{\rhohat3}^{*}\pitilde(d)_{1\rho}V_{\rhohat3}\right).\label{1709}
\end{align}
\end{Prop}
\begin{pf}
For any $\omega'\in A_{*}$ we have
\[
\begin{array}{r@{\;=\;}l}
(\id{K}\tens\,\omega'\tens\,\omega)(\id{K}\tens\,\Delta)\pitilde(d)&
(\id{K}\tens\,\omega'*\omega)\pitilde(d)=\pi(\alpha(\omega'*\omega,d))\\\Vs{5}&\pi(\alpha(\omega',\alpha(\omega,d)))=(\id{K}\tens\,\omega')\pitilde(\alpha(\omega,d)).
\end{array}
\]
Therefore
\[
\etyk{1626}
\pitilde(\alpha(\omega,d))=
(\id{K}\tens\id{A}\tens\,\omega)(\id{K}\tens\,\Delta)\pitilde(d).
\]
First formula of \rf{impleme} shows that
\[
(\id{K}\tens\,\sigma\tens\id{A})(\id{K}\tens\,\Delta)\pitilde(d)=
V_{\sigmahat3}\pitilde(d)_{1\sigma}V\Vs{3}_{\sigmahat3}^{*}.
\]
Comparing this formula with \rf{1626} we get \rf{1708}.\vs

Third formula of \rf{impleme} shows that
\[
(\id{A}\tens\,\mu\comp\rho)\Delta(a)
=(\mu\tens\id{A})\left(V\Vs{3}_{\rhohat2}^{*}\left(\rho(a)\tens\I{A}\right)V_{\rhohat2}\right)
\]
for any $a\in A$. Therefore
\[
\left(\id{K}\tens\,(\id{A}\tens\,\mu\comp\rho)\Delta\right)\pitilde(d)
=(\id{K}\tens\,\mu\tens\id{A})\left(V\Vs{3}_{\rhohat3}^{*}\pitilde(d)_{1\rho}V_{\rhohat3}\right).
\]
On the other hand, inserting $\omega=\mu\comp\rho$ in \rf{1626} we obtain
\[
\pitilde(\alpha(\mu\comp\rho,d))=
(\id{K}\tens\,(\id{A}\tens\,\mu\comp\rho)\Delta)\pitilde(d).
\]
and \rf{1709} follows.
\end{pf}

\begin{Prop}
\label{HPresults}
Let $(D,\alpha)$ be a weak $G$-dynamical system and $D_{1}\subset D$ be the range of $\alpha$ {\rm (cf \rf{910})}. Morover let $(\pi,\pitilde)$ be a faithful presentation of $\alpha$ on a Hilbert space $K$ and $(\sigma,\sigmahat)$ be a Heisenberg pair acting on a Hilbert space $H$. Then\vs

1.\ $\left(\I{K}\tens\,\sigmahat(\Ahat)\right)\pitilde(D)_{1\sigma}$ is a $\cst$-algebra of operators acting on $K\tens H$ in a non-degenerate way, \vs

2.\ We have
\[
\etyk{1812}
\left(\I{K}\tens\,\sigmahat(\Ahat)\right)\pitilde(D_{1})_{1\sigma}=
\left(\I{K}\tens\,\sigmahat(\Ahat)\right)\pitilde(D)_{1\sigma}.
\]
\end{Prop}
\begin{pf}
We shall use formula \rf{1708}. Taking the closed linear span over all $\omega\in A_{*}$ and $d\in D$ we get
\[
\pitilde(D_{1})_{1\sigma}=
\Vset{(\id{K\tens H}\tens\,\omega)
\left(V_{\sigmahat3}\pitilde(D)_{1\sigma}V\Vs{3}_{\sigmahat3}^{*}\right)}{\omega\in A_{*}}.
\]
We compute left hand side of \rf{1812}. We know that unitary $V$ belongs to $M(\Ahat\tens A)$. Therefore $(\Ahat\tens A)V=\Ahat\tens A$. In the following computation we use {\it emit-absorb} rules introduced by \rf{facto}. At first (in the second equality) $\omega$ emits $A$ to the right, later (in the forth equality) $\omega$ absorbs $A$ back. Finally in the last equality follows from second formula of \rf{algabry}.
\[
\begin{array}{r@{\;=\;}l}
\left(\I{K}\tens\,\sigmahat(\Ahat)\right)\pitilde(D_{1})_{1\sigma}
&
\Vset{(\id{K\tens H}\tens\,\omega)\left(
\left(\I{K}\tens\,\sigmahat(\Ahat)\tens\I{A}\right)
V_{\sigmahat3}\pitilde(D)_{1\sigma}V\Vs{3}_{\sigmahat3}^{*}
\right)}
{\omega\in A_{*}}\\
\Vs{6}&
\Vset{(\id{K\tens H}\tens\,\omega)\left(
\left((\Ahat\tens A)V\right)_{\sigmahat3}\pitilde(D)_{1\sigma}V\Vs{3}_{\sigmahat3}^{*}
\right)}
{\omega\in A_{*}}\\
\Vs{6}&
\Vset{(\id{K\tens H}\tens\,\omega)\left(
(\Ahat\tens A)_{\sigmahat3}\pitilde(D)_{1\sigma}V\Vs{3}_{\sigmahat3}^{*}
\right)}
{\omega\in A_{*}}\\
\Vs{6}&
\left(\I{K}\tens\,\sigmahat(\Ahat)\right)\pitilde(D)_{1\sigma}\Vset{(\id{K\tens H}\tens\,\omega)
\left(V\Vs{3}_{\sigmahat3}^{*}\right)}
{\omega\in A_{*}}\\
\Vs{6}&
\left(\I{K}\tens\,\sigmahat(\Ahat)\right)\pitilde(D)_{1\sigma}\left(\I{K}\tens\,\sigmahat(\Ahat)\right).
\end{array}
\]
We showed that
\[
\etyk{1810}
\left(\I{K}\tens\,\sigmahat(\Ahat)\right)\pitilde(D_{1})_{1\sigma}=
\left(\I{K}\tens\,\sigmahat(\Ahat)\right)\pitilde(D)_{1\sigma}\left(\I{K}\tens\,\sigmahat(\Ahat)\right).
\]
The expression on the right hand side is invariant under hermitian conjugation. Therefore
\[
\left(\I{K}\tens\,\sigmahat(\Ahat)\right)\pitilde(D_{1})_{1\sigma}=
\pitilde(D_{1})_{1\sigma}\left(\I{K}\tens\,\sigmahat(\Ahat)\right).
\]
Iterating this formula we obtain
\[
\left(\I{K}\tens\,\sigmahat(\Ahat)\right)\pitilde(D_{n})_{1\sigma}=
\pitilde(D_{n})_{1\sigma}\left(\I{K}\tens\,\sigmahat(\Ahat)\right),
\]
where $n=1,2,3\dots$ and $D_{n}$ is the product of $n$ copies of $D_{1}$. Taking closed linear span over all $n$ and remembering that $D$ is generated by $D_{1}$ we see that
\[
\etyk{1811}
\left(\I{K}\tens\,\sigmahat(\Ahat)\right)\pitilde(D)_{1\sigma}=
\pitilde(D)_{1\sigma}\left(\I{K}\tens\,\sigmahat(\Ahat)\right).
\]
Proposition \ref{XYYX} shows that 
\[
\left(\I{K}\tens\,\sigmahat(\Ahat)\right)\pitilde(D)_{1\sigma}\in\cst(K\tens H).
\]
Combining \rf{1811} with \rf{1810} we obtain \rf{1812}.
\end{pf}\vs

%The $\cst$-algebra $(\I{K}\tens\,\sigmahat(\Ahat))\pitilde(D)_{1\sigma}$ is called the crossed product of $D$ by the weak action $\alpha$ and denoted by $D\rtimes_{\alpha}\Ghat$. It turns out that it is independent of the choice of Heisenberg pair $(\sigma,\sigmahat)$ and faithful presentation $(\pi,\pitilde)$. In the following we use the notation introduced in Proposition \ref{HPresults}. We know that $\cst$-algebras $\I{K}\tens\,\sigmahat(\Ahat)$ and $\pitilde(D)_{1\sigma}$ are contained in $\M(D\rtimes_{\alpha}\Ghat)$. \vs 
%
%\begin{Prop}
%\label{nik}
%Let $(\sigma',\sigmahat')$ be another Heisenberg pair acting on a Hilbert space $H'$ and $(\pi',\pitilde')$ be another faithful presentation of $\alpha$ on a Hilbert space $K'$. Then there exists unique $\cst$-algebra iso\-morphism
%\[
%\jmath:\left(\I{K}\tens\,\sigmahat(\Ahat)\right)\pitilde(D)_{1\sigma}
%\longrightarrow
%\left(\I{K'}\tens\,\sigmahat'(\Ahat)\right)\pitilde'(D)_{1\sigma'}
%\]
%such that
%\[
%\begin{array}{r@{\;=\;}l}
%\jmath\left(\I{K}\tens\,\sigmahat(\ahat)\right)
%&
%\I{K'}\tens\,\sigmahat'(\ahat),\\
%\Vs{5}\jmath:(\pitilde(d)_{1\sigma})
%&
%\pitilde'(d)_{1\sigma'}
%\end{array}
%\]
%for any $\ahat\in\Ahat$ and $d\in D$.
%\end{Prop}

\begin{Prop}
\label{aHPresults}
Let $(D,\alpha)$ be a weak $G$-dynamical system, $D_{1}\subset D$ be the range of $\alpha$ {\rm (cf \rf{910})} and $(\pi,\pitilde)$ be a faithful presentation of $\alpha$ on a Hilbert space $K$. Then
\[
\etyk{Te1}
(\I{K}\tens A)\pitilde(D_{1})(\I{K}\tens A)=\pi(D_{1})\tens A.
\]
Moreover, if $G$ is a regular group {\rm (cf. \rf{reg})} then
\[
\etyk{Te2}
\pitilde(D_{1})(\I{K}\tens A)=\pi(D_{1})\tens A.
\]
\end{Prop}
\begin{pf}
Let $(\rho,\rhohat)$ be an anti-Heisenberg pair acting on a Hilbert space $\Hbar$. We shall use formula \rf{1709}. Taking  the closed linear span over all $\mu\in\B(\Hbar)_{*}$ and $d\in D$ we obtain
\[
\etyk{1419}
\pitilde(D_{1})=
\Vset{(\id{K}\tens\,\mu\tens\id{A})
\left(V\Vs{3}_{\rhohat3}^{*}\pitilde(D)_{1\rho}V_{\rhohat3}\right)}{\mu\in\B(\Hbar)_{*}}.
\]

 We compute $\pitilde(D_{1})(\I{K}\tens A)$. We know that unitary $V$ belongs to $M(\Ahat\tens A)$. Therefore $V(\Ahat\tens A)=\Ahat\tens A$. In the following computation we use {\it emit-absorb} rules introduced by \rf{facto}. At first (in the second equality) $\mu$ emits $\rhohat(\Ahat)$ to the left, later (in the forth equality) $\mu$ absorbs $\rhohat(\Ahat)$ back. Finally in the last equality $\mu$ emits $\rhohat(\Ahat)$ to the right.
\[
\begin{array}{r@{\;=\;}l}
\pitilde(D_{1})(\I{K}\tens A)&
\Vset{(\id{K}\tens\,\mu\tens\id{A})
\left(V\Vs{3}_{\rhohat3}^{*}\pitilde(D)_{1\rho}V_{\rhohat3}(\I{K\tens\Hbar}\tens A)\right)}{\mu\in\B(\Hbar)_{*}}\\
\Vs{7}&
\Vset{(\id{K}\tens\,\mu\tens\id{A})
\left(V\Vs{3}_{\rhohat3}^{*}\pitilde(D)_{1\rho}\left(V(\Ahat\tens A)\right)_{\rhohat3}\right)}{\mu\in\B(\Hbar)_{*}}\\
\Vs{7}&
\Vset{(\id{K}\tens\,\mu\tens\id{A})
\left(V\Vs{3}_{\rhohat3}^{*}\pitilde(D)_{1\rho}(\Ahat\tens A)_{\rhohat3}\right)}{\mu\in\B(\Hbar)_{*}}\\
\Vs{7}&
\Vset{(\id{K}\tens\,\mu\tens\id{A})
\left(\left(V\Vs{3}^{*}(\I{\Ahat}\tens A)\right)_{\rhohat3}\pitilde(D)_{1\rho}\right)}{\mu\in\B(\Hbar)_{*}}\\
\Vs{7}&
\Vset{(\id{K}\tens\,\mu\tens\id{A})
\left(\left((\Ahat\tens\I{A})V\Vs{3}^{*}(\I{\Ahat}\tens A)\right)_{\rhohat3}\pitilde(D)_{1\rho}\right)}{\mu\in\B(\Hbar)_{*}}.
\end{array}
\]
Let $X$ be a subset of $\M(\Ahat\tens A)$ introduced by formula
\[
X=(\Ahat\tens\I{A})V\Vs{3}^{*}(\I{\Ahat}\tens A).
\]
Then
\[
\etyk{oszu}
\pitilde(D_{1})(\I{K}\tens A)=\Vset{(\id{K}\tens\,\mu\tens\id{A})
\left(\Vs{3.5}X_{\rhohat3}\pitilde(D)_{1\rho}\right)}{\mu\in\B(\Hbar)_{*}}.
\]
Denote by $\Phi[X]$ the right hand side of \rf{oszu}. It makes sense for any closed linear subspace $X$ of $\M(\Ahat\tens A)$. In particular
\[
\etyk{Te3}
\begin{array}{r@{\;=\;}l}
\Phi[\Ahat\tens A]
&
\Vset{(\id{K}\tens\,\mu\tens\id{A})
\left((\Ahat\tens A)_{\rhohat3}\pitilde(D)_{1\rho}\right)}{\mu\in\B(\Hbar)_{*}}\\
\Vs{7}&\Vset{(\id{K}\tens\,\mu\tens\id{A})
\left(\Vs{3.5}(\I{K\tens\Hbar}\tens A)\pitilde(D)_{1\rho}\right)}{\mu\in\B(\Hbar)_{*}}
\\
\Vs{7}&\Vset{(\id{K}\tens\,\mu\comp\rho)
\pitilde(D)}{\Vs{4}\mu\in\B(\Hbar)_{*}}\tens A=\pi(D_{1})\tens A.
\end{array}
\]
In the second equality $\mu$ absorbs $\rhohat(\Ahat)$ from the right, at the end we use \rf{1248}. If $G$ is regular then $X=\Ahat\tens A$ (cf \rf{reg1}) and \rf{Te2} follows immediately from \rf{Te3}.\vs

We go back to general case. To compute left hand side of \rf{Te1} we use \rf{oszu}:
\[
\begin{array}{r@{\;=\;}l}
(\I{K}\tens A)\pitilde(D_{1})(\I{K}\tens A)&
\Vset{(\id{K}\tens\,\mu\tens\id{A})
\left(\Vs{3.5}\left((\I{\Ahat}\tens A)X\right)_{\rhohat3}\pitilde(D)_{1\rho}\right)}{\mu\in\B(\Hbar)_{*}}\\
\Vs{5}&\Phi[(\I{\Ahat}\tens A)X].
\end{array}
\]
On the other hand $(\I{\Ahat}\tens A)X=(\Ahat\tens A)V^{*}(\I{\Ahat}\tens A)=\Ahat\tens A$. Using again \rf{Te3} we get \rf{Te1}. 
\end{pf}\vs

Now we are able to present the prove of Baaj-Skandalis-Vaes theorem.

\begin{pf}[ of Theorem \ref{BSV}]
Let $(\pi,\pitilde)$ be a faithful presentation of $\alpha$ on a Hilbert space $K$. Identifying $D$ with $\pi(D)$ we may assume that $D\in\cst(K)$ and that $\pi$ is the embedding $D\hookrightarrow\B(K)$. In this context we write $\alphatilde$ instead of $\pitilde$. Then \rf{Te2} takes the form
\[
\alphatilde(D_{1})(\I{K}\tens A)=D_{1}\tens A.
\]
Let $n\in\natu$ and $D_{n}$ be the product of $n$ copies of $D_{1}$. We have: 
\[
\alphatilde(D_{2})(\I{K}\tens  A)=\alphatilde(D_{1})\alphatilde(D_{1})(\I{K}\tens  A)=\alphatilde(D_{1})(\I{K}\tens  A)(D_{1}\tens\I{K})=D_{2}\tens  A
\]
and similarly $\alphatilde(D_{n})(\I{K}\tens  A)=D_{n}\tens  A$ for any $n$. Taking the closed linear span over all natural $n$ we obtain Podleś condition:
\[
\alphatilde(D)(\I{K}\tens  A)=D\tens  A.
\]
It shows that $\alphatilde(D)\subset\M(D\tens A)$  and $\alphatilde\in\Mor(D,D\tens A)$. To end the proof it is sufficient to notice that in present notation \rf{0855} takes the form
\[
\alpha(\omega,d)=(\id{D}\tens\,\omega)\alphatilde(d).
\]
It means that action $\alpha$ is continuous.
\end{pf}\vs

%Let $\alpha$ be a weak action of $G$ on a $\cst$-algebra $D$. Assume that there exixts a $^{*}$-algebra homomorphism $\alphatilde\in\Mor(D,D\tens A)$ such that 
%\[
%\etyk{A01}
%\alpha(\omega,d)=(\id{D}\tens\,\omega)\alphatilde(d).
%\]
%for any $d\in D$ and $\omega\in A_{*}$. Then for any Hilbert space $H$ and any faithful $\pi\in\Rep(D,H)$, composition $\pitilde=(\pi\tens\id{A})\alphatilde\in\Mor(D,\B_{0}(H)\tens A)$ and $(\pi,\pitilde)$ is a faithful presentation of $\alpha$. Conversely we have
%\begin{Prop}
%\label{}
%Let $\alpha$ be a weak action of $G$ on a $\cst$-algebra $D$. Assume that for any Hilbert space $H$ and any faithful $\pi\in\Rep(D,H)$ there exixts $\pitilde\in\Mor(D,\B_{0}(H)\tens A)$ such that $(\pi,\pitilde)$ is a faithful presentation of $\alpha$. Then $\alpha$ is of the form \rf{A01}, where $\alphatilde\in\Mor(D,D\tens A)$.
%\end{Prop}
%
%The proof failed.

%\newpage

\section{$G$-products}
\label{sek5}
We shall analyse a structure of a $G$-product $(B,\beta,\psi)$. The section contains the proofs of Theorems \ref{Lan} and \ref{uniqueness}. At the beginning we follow the path elaborated by S. Vaes in \cite[Proof of Thm 6.7]{Vaes1}. Let $(\sigma,\sigmahat)$ be a Heisenberg pair acting on a Hilbert space $H$. The pair will be used in many constructions presented in this section. It is important to realise, which of them are independent of the choice of $(\sigma,\sigmahat)$. We shall discuss this issue later.
\begin{Prop}
\label{953a}
Let $(B,\beta,\psi)$ be a $G$-product, $\beta'=\flip\comp\beta$ and $\varphi\in\Mor(B,B\tens\B_{0}(H))$ be an injective morphism such that
\[
\etyk{859}
\varphi(b)=V_{\psi\sigma}\beta'(b)_{1\sigmahat}V_{\psi\sigma}^{*}
\]
for any $b\in B$. Then
\begin{align}
\varphi\comp\psi(\ahat)=&\I{B}\tens\,\sigmahat(\ahat),\label{1731f}\\
(\beta\tens\id{H})\varphi(b)=&\I{\Ahat}\tens\,\varphi(b)\label{1731s}
\end{align}
for any $\ahat\in\Ahat$ and $b\in B$.
\end{Prop}
\begin{pf}
Let $\ahat\in\Ahat$. According to \rf{a1}, $\beta(\psi(\ahat))=\Deltahat(\ahat)_{1\psi}$. Therefore $\beta'(\psi(\ahat))=(\flip\comp\Deltahat(\ahat))_{\psi2}$ and using the last equality in \rf{impleme} we get
\[
\beta'(\psi(\ahat))_{1\sigmahat}=(\flip\comp\Deltahat(\ahat))_{\psi\sigmahat}
=V\Vs{3}_{\psi\sigma}^{*}\left(\I{B}\tens\,\sigmahat(\ahat)\right)V\Vs{3}_{\psi\sigma}.
\]
Now \rf{1731f} follows immediately from \rf{859}.\vs

Let $b\in B$. Then
\[
\etyk{dorazny}
(\beta\tens\id{H})\varphi(b)=(\beta\tens\id{H})\left(V_{\psi\sigma}\right)\;(\beta\tens\id{H})\left(\beta'(b)\Vs{3.5}_{1\sigmahat}\right)\;(\beta\tens\id{H})\left(V_{\psi\sigma}^{*}\right).
\]
Taking into account \rf{a1} and second formula of \rf{1037} we obtain:
\[
\begin{array}{r@{\;=\;}l}
(\beta\tens\id{H})\left(V_{\psi\sigma}\Vs{3.5}\right)& (\beta\comp\psi\tens\,\sigma)V=
\left((\id{\Ahat}\tens\,\psi)\Deltahat\tens\,\sigma\right)V\\
\Vs{5}&\left(\id{\Ahat}\tens\,\psi\tens\,\sigma\right)\left(\Deltahat\tens\id{A}\right)V=
\left(\id{\Ahat}\tens\,\psi\tens\,\sigma\right)\left(V_{23}V_{13}\right)=
V_{\psi\sigma}V_{1\sigma}.
\end{array}
 \]
 We know that $\beta$ is a left action of $\Ghat$: $(\id{\Ahat}\tens\beta)\beta=(\Deltahat\tens\id{B})\beta$. A moment of reflection shows that 
 $(\beta\tens\id{\Ahat})\beta'=\flip_{12}(\id{B}\tens\flip\Deltahat)\beta'$. Taking into account the last formula of \rf{impleme} we obtain
\[
\begin{array}{r@{\;=\;}l}
(\beta\tens\id{H})\left(\beta'(b)\Vs{3.5}_{1\sigmahat}\right)&
(\beta\tens\,\sigmahat)\beta'(b)=
(\id{\B_{0}(H)\tens B}\tens\,\sigmahat)\flip_{12}(\id{B}\tens\flip\Deltahat)\beta'(b)\\
&
\flip_{12}(\id{B}\tens\,(\id{\Ahat}\tens\,\sigmahat)\flip\Deltahat)\beta'(b)\Vs{5}\\
\Vs{6}&
\flip_{12}\left(V\Vs{3}_{2\sigma}^{*}\beta'(b)_{1\sigmahat}V_{2\sigma}\right)
=V_{1\sigma}^{*}\beta'(b)_{2\sigmahat}V_{1\sigma}.
\end{array}
 \]
Inserting the results of above computations to \rf{dorazny} we see that
 \[
 \begin{array}{r@{\;=\;}l}
(\beta\tens\id{H})\varphi(b)&V_{\psi\sigma}V_{1\sigma}V_{1\sigma}^{*}\beta'(b)_{2\sigmahat}V_{1\sigma}\;V_{1\sigma}^{*}V_{\psi\sigma}^{*}\\
 &V_{\psi\sigma}\beta'(b)_{2\sigmahat}V_{\psi\sigma}^{*}=\varphi(b)_{23}\Vs{5}=
 \I{B}\tens\,\varphi(b).
\end{array}
 \]
Formula \rf{1731f} is proven
\end{pf}
\begin{Prop}
\label{953b}
Let $(B,\beta,\psi)$ be a $G$-product, $\varphi\in\Mor(B,B\tens\B_{0}(H))$ be the injective morphism introduced by \rf{859} and
\[
\etyk{937}
D_{1}=\Vset{(\id{B}\tens\,\mu)\varphi(B)}{\Vs{4}\mu\in\B(H)_{*}}.
\]
Then $D_{1}$ is a norm closed linear subset of $\M(B)$. Elements of $D_{1}$ are $\beta$-invariant: 
\[
\etyk{3441}
\beta(d)=\I{\Ahat}\tens\,d
\]
 for any $d\in D_{1}$ and
\[
\etyk{3442}
B=\psi(\Ahat)D_{1}.
\]
\end{Prop}

\begin{pf} Relation \rf{3441} follows immediately from \rf{1731s}.\vs

We know that $V_{1\sigma}$ is a unitary multiplier of $\Ahat\tens\B_{0}(H)$. Therefore $V_{\psi\sigma}$ is a unitary multiplier of $B\tens\B_{0}(H)$ and
\[
\begin{array}{r@{\;=\;}l}
(\Ahat\tens\B_{0}(H))V_{1\sigma}&\Ahat\tens\B_{0}(H)\\
\Vs{4}(B\tens\B_{0}(H))V_{\psi\sigma}^{*}&B\tens\B_{0}(H)
\end{array}
\]
Moreover Podleś condition \rf{1733} easily shows that
\[
(\I{B}\tens\B_{0}(H))\beta'(B)_{1\sigmahat}=B\tens\B_{0}(H).
\]
Using the above formulae we obtain
 \[
 \begin{array}{r@{\;=\;}l}
(\psi(\Ahat)\tens\B_{0}(H))\varphi(B)
&(\psi(\Ahat)\tens\B_{0}(H))V_{\psi\sigma}\beta'(B)_{1\sigmahat}V_{\psi\sigma}^{*}\\
\Vs{5}&(\psi\tens\id{H})\left((\Ahat\tens\B_{0}(H))V_{1\sigma}\right)\beta'(B)_{1\sigmahat}V_{\psi\sigma}^{*}
 \\
\Vs{5}&(\psi\tens\id{H})\left(\Ahat\tens\B_{0}(H)\right)\beta'(B)_{1\sigmahat}V_{\psi\sigma}^{*}
 \\
\Vs{5}&(\psi(\Ahat)\tens\I{H})\left(\I{B}\tens\B_{0}(H)\right)\beta'(B)_{1\sigmahat}V_{\psi\sigma}^{*}
 \\
 \Vs{5}&(\psi(\Ahat)\tens\I{H})\left(B\tens\B_{0}(H)\right)V_{\psi\sigma}^{*}
=(\psi(\Ahat)\tens\I{H})\left(B\tens\B_{0}(H)\right)
 \\
\Vs{5}&\psi(\Ahat)B\tens\B_{0}(H).
\end{array}
 \]
 Taking into account the equality $\psi(\Ahat)B=B$ (this is because $\psi\in\Mor(\Ahat,B)$) we get
 \[
(\psi(\Ahat)\tens\B_{0}(H))\varphi(B)=B\tens\B_{0}(H).
\]
Now, using  \rf{facto} we obtain
\[
 \begin{array}{r@{\;=\;}l}
 \psi(\Ahat)D_{1}&\Vset{\psi(\Ahat)(\id{B}\tens\,\omega)\varphi(B)}{\Vs{4}\omega\in\B(H)_{*}}\\
 \Vs{8}&\Vset{(\id{B}\tens\,\omega)\left\{(\psi(\Ahat)\tens\B_{0}(H)\Vs{3.5})\varphi(B)\right\}}{\Vs{4}\omega\in\B(H)_{*}}\\
  \Vs{8}&\Vset{(\id{B}\tens\,\omega)\left\{B\tens\B_{0}(H)\Vs{3.5}\right\}}{\Vs{4}\omega\in\B(H)_{*}}=B.
\end{array}
 \]
Formula \rf{3442} is proven. 
\end{pf}\vs

With the notation introduced in the previous Propositions we have
\begin{Prop}
\label{LAEU}
$\cst$-algebra generated by $D_{1}$ is the only Landstad algebra of $(B,\beta,\psi)$.
\end{Prop}
\begin{pf}
Let $D$ be a $\cst$-subalgebra of $\M(B)$ satisfying Conditions C1 and C2 of Definition \ref{LC}. Condition C1 shows that $\beta(D)=\I{\Ahat}\tens D$. Therefore $\varphi(D)=(\id{B}\tens\,\sigma)\left(V_{\psi2}(D\tens\I{A})V_{\psi2}^{*}\right)$. $\I{\Ahat}$ was replaced by $\I{A}$ because $\sigmahat(\I{\Ahat})=\I{H}=\sigma(\I{A})$. We know that $\omega\in A_{*}$ if and only if $\omega=\mu\comp\sigma$, where $\mu\in\B(H)_{*}$. Therefore
\[
\Vset{(\id{B}\tens\,\omega)\left(\Vs{3.5}V_{\psi2}(D\tens\I{\Ahat})V_{\psi2}^{*}\right)}{\omega\in A_{*}}
=\Vset{(\id{B}\tens\,\mu)\varphi(D)}{\Vs{4}\mu\in\B(H)_{*}}
\]
We shall prove that the latter set coincides with $D_{1}$. Indeed:
\[
\begin{array}{r@{\;=\;}l}
\Vset{(\id{B}\tens\,\mu)\varphi(D)}{\Vs{4}\mu\in\B(H)_{*}}
&\Vset{(\id{B}\tens\,\mu)
\left(\left(\I{B}\tens\,\sigmahat(\Ahat)\right)\varphi(D)\right)}{\Vs{2}\mu\in\B(H)_{*}}%
\\
\Vs{8}&\Vset{(\id{B}\tens\,\mu)\varphi\left(\psi(\Ahat)\Vs{3.5}D\right)}{\Vs{4}\mu\in\B(H)_{*}}
\\
\Vs{8}&\Vset{(\id{B}\tens\,\mu)\varphi(B)}{\Vs{4}\mu\in\B(H)_{*}}=D_{1}.
\end{array}
\]
The first equality follows from \rf{facto}, the second one from \rf{1731f} and the third one from C2.
We showed that 
\[
\etyk{0555}
\Vset{(\id{B}\tens\,\omega)\left(\Vs{3.5}V_{\psi2}(D\tens\I{\Ahat})V_{\psi2}^{*}\right)}{\omega\in A_{*}}
=D_{1}.
\]

Assume now that $D$ is a Landstad algebra for $(B,\beta,\psi)$. Then Condition C3 shows that $D$ is generated by $D_{1}$.\vs

Conversely assume that $D$ is the $\cst$-algebra generated by $D_{1}$. Then $D\subset\M(B)$. A moment of reflection shows that
\[
D=\left(D_{1}\cup\Vs{4}D_{1}D\right)^{\rm CLS}.
\]
We know that $\psi(\Ahat)D_{1}=B$ (see \rf{3442}). The simple observation: $\psi(\Ahat)D_{1}D=BD\subset B$ combined with the above formula shows that $\psi(\Ahat)D=B$. Condition C2 is verified.\vs

According to \rf{3441} elements of $D_{1}$ are $\beta$-invariant. The same is true for elements of $D$ and Condition C1 is verified. \vs

Now we may use \rf{0555}. Remembering that $D$ is generated by $D_{1}$ we see that Condition C3 is satisfied. It means that $D$ is a Landstad algebra for $(B,\beta,\psi)$.\vs

The existence and uniqueness of Landstad algebra is shown.
\end{pf}
\vs

Now we are able to prove Theorem \ref{Lan}.\vs

\begin{pf} Let $D\subset\M(B)$ be the Landstad algebra of $G$-product $(B,\beta,\psi)$ (see Definition \ref{LC}). It follows immediately from C3, that formula \rf{1903} defines a mapping $\alpha$ acting from $ A_{*}\times D$ into $D$ satisfying Condition WA1 of Definition \ref{mea}. We shall check Condition WA0.\vs

One  may assume that $B\in\cst(K)$, where $K$ is a Hilbert space. Then $D\subset\M(B)\subset\B(K)$. 
Denote by $\pi'$ the embedding of $B$ into $\B(K)$, by $\pi$ the restriction of $\pi'$ to $D$ (so $\pi$ is the embedding of $D$ into $\B(K)$) and by $\pitilde\in\Mor(D,\B_{0}(K)\tens A)$ the morphism introduced by the formula
\[
\pitilde(d)=(\pi'\tens\id{H})\left(V_{\psi2}(d\tens\I{H})V_{\psi2}^{*}\right)
\]
Let $\omega\in A_{*}$. Taking into account \rf{1903} we obtain
\[
\begin{array}{r@{\;=\;}l}
\pi(\alpha(\omega,d))&(\pi'\tens\,\omega)\left(V_{\psi2}(d\tens\I{H})V_{\psi2}^{*}\right)\\
\Vs{4}&(\id{K}\tens\,\omega)\pitilde(d).
\end{array}
\]
This formula coincides with \rf{0855}. Therefore $(\pi,\pitilde)$ is a faithful presentation of $\alpha$ and Condition WA0 is verified.\vs

We shall show that $\alpha$ introduced by \rf{1903} obeys Condition WA2. This is a matter of easy computation. Let $\omega,\omega'\in A_{*}$ and $d\in D$. Then
\[
\etyk{vvv}
\begin{array}{r@{\;=\;}l}
\alpha(\omega*\omega',d)&(\id{B}\tens\,(\omega*\omega'))\left(V_{\psi2}(d\tens\I{A})V_{\psi2}^{*}\right)\\
\Vs{5}&(\id{B}\tens\,\omega\tens\,\omega')(\id{B}\tens\,\Delta)\left(V_{\psi2}(d\tens\I{A})V_{\psi2}^{*}\right)\\
\Vs{5}&(\id{B}\tens\,\omega\tens\,\omega')\left(V_{\psi2}V_{\psi3}(d\tens\I{A}\tens\I{A})
V_{\psi3}^{*}V_{\psi2}^{*}\right)\\
\Vs{5}&(\id{B}\tens\,\omega)\left(V_{\psi2}(\alpha(\omega',d)\tens\I{A})V_{\psi2}^{*}\right)=
\alpha(\omega,\alpha(\omega',d)).
\end{array}
\]
We have shown that $\alpha$ is a weak right action of $G$ on $D$. The proof of Theorem \ref{Lan} is complete.
\end{pf}\vs

For any $d\in D$ and $\mu\in\B_{0}(H)$ we have
\[
\etyk{2056}
(\id{B}\tens\,\mu)\varphi(d)=\alpha(\mu\comp\sigma,d).
\]
Indeed, $\beta(d)=\I{\Ahat}\tens\,d$, $\beta'(d)=d\tens\I{\Ahat}$, $\beta'(d)_{1\sigmahat}=d\tens\I{H}$ and $\varphi(d)=V_{\psi\sigma}(d\tens\I{H})V
\Vs{3}_{\psi\sigma}^{*}$. Therefore 
\[
(\id{B}\tens\,\mu)\varphi(d)=(\id{B}\tens\,\mu\comp\sigma)\left(V_{\psi2}(d\tens\I{A})V\Vs{3}_{\psi2}^{*}\right)
\]
and using \rf{1903} we get \rf{2056}.
\vs

To prove the uniqueness Theorem \ref{uniqueness} we start with the following proposition:
\begin{Prop}
\label{isometry}
Let $(B^{1},\beta^{1},\psi^{1})$ and $(B^{2},\beta^{2},\psi^{2})$ be $G$-products and $(D^{1},\alpha^{1})$, $(D^{2},\alpha^{2})$ be corresponding Landstad weak $G$-dynamical systems and $\imath:D^{1}\longrightarrow D^{2}$ be an isomorphism of the weak $G$-dynamical systems $(D^{1},\alpha^{1})$ and $(D^{2},\alpha^{2})$. It means that
\[
\etyk{2718}
\imath(\alpha^{1}(\omega,d))=\alpha^{2}(\omega,\imath(d))
\]
for any $\omega\in A_{*}$ and $d\in D^{1}$. Then $B^{1}$ and $B^{2}$ are isomorphic as Banach spaces. More precisely: there exists a linear, isometric map $\jmath:B^{1}\longrightarrow B^{2}$ such that $\jmath(B^{1})=B^{2}$ and
\[
\etyk{2050}
\begin{array}{r@{\;=\;}l}
\jmath\left(\psi^{1}(\ahat)d\right)\Vs{5}&\psi^{2}(\ahat)\imath(d),\\
\jmath\left(d\psi^{1}(\ahat')\right)\Vs{5}&\imath(d)\psi^{2}(\ahat'),\\
\jmath\left(\psi^{1}(\ahat)d\psi^{1}(\ahat')\right)\Vs{5}&\psi^{2}(\ahat)\imath(d)\psi^{2}(\ahat')
\end{array}
\]
for any $\ahat,\ahat'\in\Ahat$ and $d\in D^{1}$.
\end{Prop}
\begin{pf} Let $D$ be the Landstad algebra of a $G$-product $(B,\beta,\psi)$. We know that $\psi(\Ahat)D=B$. Applying to the both sides hermitian conjugation we get $D\psi(\Ahat)=B$. Moreover $\psi(\Ahat)D\psi(\Ahat)=B$. Indeed $\psi\in\Mor(\Ahat,B)$. Therefore $\psi(\Ahat)B=B$ and $\psi(\Ahat)D\psi(\Ahat)=\psi(\Ahat)B=B$. It shows that each of the three formulas \rf{2050} defines $\jmath$ on a dense subset of $B^{1}$ with values running over a dense subset of $B^{2}$.\vs

Let $(\sigma,\sigmahat)$ be a Heisenberg pair acting on a Hillbert space $H$, $\varphi\in\Mor(B,B\tens\B_{0}(H))$  be the morphism introduced by \rf{859} and $\alpha$ be the weak action of $G$ on $D$ introduced by \rf{1903}. We shall use the $\B(H)$-bimodul structure of $\B(H)_{*}$. According to \rf{1410}, for any $\mu\in\B(H)_{*}$ and $\ahat,\ahat'\in\Ahat$, $\mu\sigmahat(\ahat),\ \sigmahat(\ahat')\mu$ and $\sigmahat(\ahat')\mu\sigmahat(\ahat)$ are normal linear functionals on $\B(H)$ such that
\[
\begin{array}{r@{\;=\;}l}
\left(\mu\sigmahat(\ahat)\right)(m)&\mu(\sigmahat(\ahat)m),\\
\Vs{5}\left(\sigmahat(\ahat')\mu\right)(m)&\mu(m\sigmahat(\ahat')),\\  
\Vs{5}\left(\sigmahat(\ahat')\mu\sigmahat(\ahat)\right)(m)&\mu(\sigmahat(\ahat)m\sigmahat(\ahat'))
\end{array}
\]
for any $m\in\B(H)$. Let $d\in D$. We shall prove that 
\[
\etyk{2717}
\begin{array}{r@{\;=\;}l}
\alpha((\mu\sigmahat(\ahat))\comp\sigma,d)
&
(\id{B}\tens\,\mu)\varphi(\psi(\ahat)d),
\\
\alpha((\sigmahat(\ahat')\mu)\comp\sigma,d)
&
\Vs{5}(\id{B}\tens\,\mu)\varphi(d\psi(\ahat'))
\\
\alpha((\sigmahat(\ahat')\mu\sigmahat(\ahat))\comp\sigma,d)
&
\Vs{5}(\id{B}\tens\,\mu)\varphi(\psi(\ahat)d\psi(\ahat')).
\end{array}
\]
Indeed taking into account \rf{2056} and \rf{1731f} we have
\[
\begin{array}{r@{\;=\;}l}
\alpha((\mu\sigmahat(\ahat))\comp\sigma,d)&
(\id{B}\tens\,\mu\sigmahat(\ahat))\varphi(d)=
(\id{B}\tens\,\mu)\left((\I{B}\tens\,\sigmahat(\ahat))\varphi(d)\Vs{3}\right)
\\ \Vs{5}&(\id{B}\tens\,\mu)\left(\varphi(\psi(\ahat))\varphi(d)\Vs{3}\right)=
(\id{B}\tens\,\mu)\varphi(\psi(\ahat)d).
\end{array}
\]
In the same way one can prove two remaining formulas of \rf{2717}.\vs

Let $b^{1}\in B^{1}$ and $b^{2}\in B^{2}$ be elements of the form
\[
\begin{array}{r@{\;=\;}c@{\;+\;}c@{\;+\;}c}
b^{1}&{\displaystyle \sum_{k=1}^{r-1}\psi^{1}(\ahat_{k})d_{k}}
&
{\displaystyle \sum_{k=r}^{s-1}d_{k}\psi^{1}(\ahat'_{k})}
&
{\displaystyle \sum_{k=s}^{t-1}\psi^{1}(\ahat_{k})d_{k}\psi^{1}(\ahat'_{k})},\\
\Vs{8}b^{2}
&
{\displaystyle \sum_{k=1}^{r-1}\psi^{2}(\ahat_{k})\imath(d_{k})}
&
{\displaystyle \sum_{k=r}^{s-1}\imath(d_{k})\psi^{2}(\ahat'_{k})}
&
{\displaystyle \sum_{k=s}^{t-1}\psi^{2}(\ahat_{k})\imath(d_{k})\psi^{2}(\ahat'_{k})},
\end{array}
\] 
\Vs{4}where $1\leq r\leq s\leq t$ are integers, $d_{k}\in D^{1}$ for $1\leq k< t$, $\ahat_{k}\in\Ahat$ for $1\leq k\leq r$ and $s\leq k<t$ and $\ahat'_{k}\in\Ahat$ for $r\leq k<t$. We look for an isometry $\jmath: B^{1}\longrightarrow B^{2}$ such that $\jmath(b^{1})=b^{2}$. To prove the existence of $\jmath$ it is sufficient to show that
\[
\etyk{nn}
\norm{b^{1}}=\norm{b^{2}}.
\]

Let $\varphi^{1}\in\Mor(B^{1},B^{1}\tens\B_{0}(H))$ and $\varphi^{2}\in\Mor(B^{2},B^{2}\tens\B_{0}(H))$ be morhisms  related to $G$-products $(B^{1},\beta^{1},\psi^{1})$ and $(B^{2},\beta^{2},\psi^{2})$ by the formula \rf{859}. Combining \rf{2717} with \rf{2718} we see that
\[
\etyk{1859}
\begin{array}{r@{\;=\;}l}
\imath\left((\id{B^{1}}\tens\,\mu)\varphi^{1}(\psi^{1}(\ahat)d)\right)
&(\id{B^{2}}\tens\,\mu)\varphi^{2}(\psi^{2}(\ahat)\imath(d))
\\
\imath\left((\id{B^{1}}\tens\,\mu)\varphi^{1}(d\psi^{1}(\ahat'))\right)
\Vs{5}&(\id{B^{2}}\tens\,\mu)\varphi^{2}(\imath(d)\psi^{2}(\ahat'))
\\
\imath\left((\id{B^{1}}\tens\,\mu)\varphi^{1}(\psi^{1}(\ahat)d\psi^{1}(\ahat'))\right)
\Vs{5}&(\id{B^{2}}\tens\,\mu)\varphi^{2}(\psi^{2}(\ahat)\imath(d)\psi^{2}(\ahat'))
\end{array}
\]
for any $\ahat,\ahat'\in\Ahat$, $d\in D^{1}$ and $\mu\in\B(H)_{*}$. It shows that
\[
\etyk{1860}
\imath\left((\id{B^{1}}\tens\,\mu)\varphi^{1}(b^{1})\right)
=(\id{B^{2}}\tens\,\mu)\varphi^{2}(b^{2})
\]
for any $\mu\in\B(H)_{*}$.\vs

Let $m,n$ be one-dimensional operators acting on $H$: $m=\ket{x}\bra{y}$ and $n=\ket{z}\bra{u}$, where $x,y,z,u$ are elements of $H$ and let $\mu$ be a normal linear functional on $\B(H)$ such that $\mu(a)=\ITS{y}{a}{z}$ for any $a\in\B(H)$.. Then
\[
\begin{array}{r@{\;=\;}l}
(\I{B^{1}}\tens\,m)\varphi^{1}(b^{1})(\I{B^{1}}\tens\,n)&(\id{B^{1}}\tens\,\mu)\varphi^{1}(b^{1})\tens\ket{x}\bra{u}\in D^{1}\tens\B_{0}(H),\\
(\I{B^{2}}\tens\,m)\varphi^{2}(b^{2})(\I{B^{2}}\tens\,n)&(\id{B^{2}}\tens\,\mu)\varphi^{2}(b^{2})\tens\ket{x}\bra{u}\in D^{2}\tens\B_{0}(H).\Vs{5}
\end{array}
\]
Formula \rf{1860} shows now that
\[
(\imath\tens\id{H})\left((\I{B^{1}}\tens\,m)\varphi^{1}(b^{1})(\I{B^{1}}\tens\,n)\right)=(\I{B^{2}}\tens\,m)\varphi^{2}(b^{2})(\I{B^{2}}\tens\,n).
\]
Any compact operator is a norm limit of sums of one-dimensional operators. Therefore the above formula holds for all $m,n\in\B_{0}(H)$. $\imath$ is an isomorphism of $\cst$-algebras. Therefore
\[
\norm{(\I{B^{1}}\tens\,m)\varphi^{1}(b^{1})(\I{B^{1}}\tens\,n)}=\norm{(\I{B^{2}}\tens\,m)\varphi^{2}(b^{2})(\I{B^{2}}\tens\,n)}.
\]
Replacing $m$ and $n$ by approximate units of $\B_{0}(H)$ and passing to the limit we get
\[
\norm{\varphi^{1}(b^{1})}=\norm{\varphi^{2}(b^{2})}.
\]
Remembering that $\varphi^{1}$ and $\varphi^{2}$ are injective morphisms we obtain \rf{nn}.
\end{pf}\vs

Now the proof of the uniqueness theorem is easy.

\begin{pf}[ of Theorem \ref{uniqueness}]
At first we show that the isometry $\jmath:B^{1}\longrightarrow B^{2}$ introduced by \rf{2050} is a $\cst$-algebra isomorphism i.e.
\[
\etyk{cstmorf}
\begin{array}{r@{\;=\;}l}
\jmath(b^{*})&\jmath(b)^{*},\\
\Vs{5}\jmath(bb')&\jmath(b)\jmath(b')
\end{array}
\]
for any $b,b'\in B^{1}$.
\vs

Let $\ahat\in\Ahat$ and $d\in D^{1}$. Then using the second and the first formula of \rf{2050} we obtain
\[
\jmath((\psi^{1}(\ahat)d)\Vs{3}^{*})=\jmath(d\Vs{3}^{*}\psi^{1}(\ahat\Vs{3}^{*}))=\imath(d\Vs{3}^{*})\psi^{2}(\ahat\Vs{3}^{*})=\left(\psi^{2}(\ahat)\imath(d)\right)\!\Vs{3}^{*}=\jmath(\psi^{1}(\ahat)d)\Vs{3}^{*}.
\] 
Remembering that $\psi^{1}(\Ahat)D^{1}=B^{1}$ we get the first formula of \rf{cstmorf}.\vs

Let $\ahat,\ahat'\in\Ahat$ and $d,d'\in D^{1}$. Then using \rf{2050} we obtain
\[
\jmath((\psi^{1}(\ahat)dd'\psi^{1}(\ahat'))=
\psi^{2}(\ahat)\imath(dd')\psi^{2}(\ahat')=
\psi^{2}(\ahat)\imath(d)\imath(d')\psi^{2}(\ahat')=
\jmath(\psi^{1}(\ahat)d)\jmath(d'\psi^{1}(\ahat')).
\] 
Remembering that $\psi^{1}(\Ahat)D^{1}=D^{1}\psi^{1}(\Ahat)=B^{1}$ we get the second formula of \rf{cstmorf}. \vs

Now we know that $\jmath\in\Mor(B^{1},B^{2})$. Extending $\jmath$ to multiplier algebras we may rewrite first formula of \rf{2050} in the following way:
\[
\jmath\left(\psi^{1}(\ahat)\right)\jmath(d)=\psi^{2}(\ahat)\imath(d).
\]
It shows that
\[
\etyk{2703}
\begin{array}{r@{\;=\;}l}
\jmath\left(\psi^{1}(\ahat)\right)&\psi^{2}(\ahat),\\
\Vs{5}\jmath(d)&\imath(d)
\end{array}
\]
for any $\ahat\in\Ahat$ and $d\in D^{1}$. Second equality says that $\jmath$ is an extension of $\imath$, the first one shows that the triangle in diagram \rf{a2} is commutative. We shall prove that the square in \rf{a2} is also commutative. Let $\ahat\in\Ahat$ and $d\in D^{!}$. Then using \rf{a1} and the triangle in \rf{a2} we get
\[
(\id{\Ahat}\tens\jmath)\beta^{1}(\psi^{1}(\ahat))=(\id{\Ahat}\tens\jmath)(\id{\Ahat}\tens\,\psi^{1})\Deltahat(\ahat)=(\id{\Ahat}\tens\,\psi^{2})\Deltahat(\ahat)= \beta^{2}(\psi^{2}(\ahat)).
\]
Similarly using Condition C1 of Definition \ref{LC} we get
\[
(\id{\Ahat}\tens\jmath)\beta^{1}(d)=(\id{\Ahat}\tens\jmath)(\I{\Ahat}\tens\,d)=
\I{\Ahat}\tens\imath(d)=\beta^{2}(\imath(d)).
\]
Combining the two formulas and using the first relation of \rf{2050} we obtain 
\[
(\id{\Ahat}\tens\jmath)\beta^{1}(\psi^{1}(\ahat)d)=\beta^{2}(\psi^{2}(\ahat)\imath(d))=\beta^{2}(\jmath(\psi^{1}(\ahat)d))
\] 
Remembering that $\psi^{1}(\Ahat)D^{1}=B^{1}$ we conclude that
\[
(\id{\Ahat}\tens\jmath)\beta^{1}(b)=\beta^{2}(\jmath(b))
\]
for any $b\in B^{1}$. It shows that \rf{a2} is a commutative diagram, so $\jmath$ is an isomorphism of $G$-products $(B^{1},\beta^{1},\psi^{1})$ and $(B^{2},\beta^{2},\psi^{2})$.\vs

Conversely if $\jmath$ is  an isomorphism of $G$-products $(B^{1},\beta^{1},\psi^{1})$ and $(B^{2},\beta^{2},\psi^{2})$ extending $\imath$ then $\jmath$ satisfies relations \rf{2703}. The latter implies the first relation of \rf{2050}. Remembering that $\psi^{1}(\Ahat)D^{1}$ coincides with $B^{1}$ we see that $\jmath$ is uniquely determined.
\end{pf}\vs

The Heisenberg pair $(\sigma,\sigmahat)$ chosen at the beginning of present section plays an important role in our considerations. Definitions of $\varphi\in\Mor(B,\B_{0}(H)\tens B)$ and of $D_{1},D\subset\M(B)$ contain explicitly $\sigma$ and $\sigmahat$. However $D$ as unique Landstad algebra does not depend on the choice of $(\sigma,\sigmahat)$. Formula \rf{0555} shows that also $D_{1}$ remains the same when we change $(\sigma,\sigmahat)$.\vs 

\section{Crossed product construction}
\label{sek4}

This section is devoted to the proof of Theorem \ref{CP}. Let $(D,\alpha)$ be a weak $G$-dynamical system, $K,H$ be Hilbert spaces, $(\sigma,\sigmahat)$ be a canonical  Heisenberg pair acting on $H$ and  $(\pi,\pitilde)$ be a faithful presentation of $\alpha$ on $K$. Then $\pitilde\in\Mor(D,\B_{0}(K)\tens A)$ is injective and 
\[
\etyk{2107}
\pi(\alpha(\omega,d))=(\id{K}\tens\,\omega)\pitilde(d)
\]
for any $\omega\in A_{*}$ and $d\in D$.\vs

For any $\ahat\in\Ahat$ we set $\psi(\ahat)=\I{K}\tens\,\sigmahat(\ahat)$. Clearly $\psi\in\Rep(\Ahat,K\tens H)$. Let 
\[
\etyk{75}
B=\psi(\Ahat)(\id{K}\tens\,\sigma)\pitilde(D).
\]
Then $B$ is a closed linear subset of $\B(K\tens H)$. We have to show that the above assumptions imply all three statements of Theorem \ref{CP}. Statement 1 is already established (see Proposition \ref{HPresults}.1).\vs

\begin{pf}[ of statement 2 of Theorem \ref{CP}]
In this section $(\sigma,\sigmahat)$ is a canonical Heisenberg pair acting on a Hilbert space $H$. Let $\rho$ be a representation of $A$ acting on $H$ such that $(\rho,\sigma)$ is a commuting pair and $(\rho,\sigmahat)$ is an anti-Heisenberg pair. We shall deal with elements of $\M(\Ahat\tens\B_{0}(K\tens H))$. Among them we have unitary $V_{1\rho}$ having the second leg (the one in $K$) trivial. We recall that $B\subset\B(K\tens H)$. For any $b\in B$ we set
\[
\beta(b)=V_{1\rho}(\I{\Ahat}\tens b)V_{1\rho}^{*}.
\]
Then $\beta$ is a $^{*}$-algebra homomorphism from $B$ to $\M(\Ahat\tens\B_{0}(K\tens H))$. We shall prove that  
\[
\etyk{1632}
\begin{array}{r@{\;=\;}l}
\beta(\psi(\ahat))&(\id{\Ahat}\tens\,\psi)\Deltahat(\ahat)\\
\beta((\id{K}\tens\,\sigma)\pitilde(d))&\I{\Ahat}\tens\,(\id{K}\tens\,\sigma)\pitilde(d)\Vs{6}
\end{array}
\]
for any $\ahat\in\Ahat$ and $d\in D$.\vs

Indeed, we know that $(\rho,\sigmahat)$ is an anti-Heisenberg pair. In this case second formula of \rf{impleme} takes the form $(\id{\Ahat}\tens\,\sigmahat)\Deltahat(\ahat)=V_{1\rho}(\I{\Ahat}\tens\,\sigmahat(\ahat))V_{1\rho}^{*}$. Reinterpreting this equation as the equality of elements of $\M(\Ahat\tens\B_{0}(K)\tens\B_{0}(H))$ having trivial second leg we obtain first formula of \rf{1632}. We know that $(\rho,\sigma)$ is a commuting pair. Therefore $V_{1\rho}$ commutes with $\I{\Ahat}\tens\,(\id{K}\tens\,\sigma)\pitilde(d)$ for all $d\in D$ and second formula of \rf{1632} follows.\vs

Combining the two formulae of \rf{1632} we get
\[
\etyk{impan0}
\beta(\psi(\ahat)\pitilde(d))=(\id{\Ahat}\tens\,\psi)\Deltahat(\ahat)(\I{\Ahat}\tens\pitilde(d))
\]
for any $\ahat\in\Ahat$ and $d\in D$. We shall prove that $\beta\in\Mor(B,\Ahat\tens B)$. Indeed  
\[
\begin{array}{r@{\;=\;}l}
(\Ahat\tens\I{B})\beta(B)
&(\Ahat\tens\I{B})(\id{\Ahat}\tens\,\psi)\Deltahat(\Ahat)(\I{\Ahat}\tens\pitilde(D))\\
\Vs{5}&(\id{\Ahat}\tens\,\psi)\left((\Ahat\tens\I{\Ahat})\Deltahat(\Ahat)\right)(\I{\Ahat}\tens\pitilde(D))\\
\Vs{5}&(\id{\Ahat}\tens\,\psi)(\Ahat\tens\Ahat)(\I{\Ahat}\tens\pitilde(D))=\Ahat\tens B.
\end{array}
\] 
So we have
\[
\etyk{Podles}
(\Ahat\tens\I{B})\beta(B)=\Ahat\tens B.
\]
Multiplying both sides (from the left) by $\I{\Ahat}\tens B$ we obtain $(\Ahat\tens B)\beta(B)=\Ahat\tens B$. It means that $\beta\in\Mor(B,A\tens B)$.\vs

First formula of \rf{1632} shows that \rf{a1} is a commutative diagram. Using commutativity of \rf{a1} and coassociativity of $\Deltahat$ one can easily show that $(\id{\Ahat}\tens\beta)\beta\comp\psi=(\Deltahat\tens\id{B})\beta\comp\psi$. Consequently 
\[
(\id{\Ahat}\tens\beta)\beta(\psi(\ahat))=(\Deltahat\tens\id{B})\beta(\psi(\ahat))
\]
for any $\ahat\in\Ahat$. Moreover for any $d\in D$ we have: 
\[
\begin{array}{r@{\;=\;}l}
(\id\tens\beta)\beta((\id{K}\tens\,\sigma)\pihat(d))&(\id\tens\beta)(\I{\Ahat}\tens\,(\id{K}\tens\,\sigma)\pihat(d))=\I{\Ahat}\tens\I{\Ahat}\tens\,(\id{K}\tens\,\sigma)\pihat(d)
\\&(\Deltahat\tens\id{B})(\I{\Ahat}\tens\,(\id{K}\tens\,\sigma)\pihat(d))=(\Deltahat\tens\id{B})\beta((\id{K}\tens\,\sigma)\pihat(d)).\Vs{5}
\end{array}
\]

 Combining two last formulae and remembering that $B=\psi(\Ahat)(\id{K}\tens\,\sigma)\pihat(D)$ we obtain that
\[
(\id{\Ahat}\tens\beta)\beta(b)=(\Deltahat\tens\id{B})\beta(b)
\]
for any $b\in B$. It means that $\beta$ is a left action of $G$ on $B$. The action is continuous: \rf{Podles} is the Podleś condition for this action. Statement 2 is shown.
\end{pf}

\begin{pf}[ of statement 3 of Theorem \ref{CP}]
We already know that the diagram \rf{a1} is commutative. Hence $(B,\beta,\psi)$ is a $G$-product.\vs

We have to show that $(\id{K}\tens\,\sigma)\pitilde(D)$ is a Landstad algebra for $(B,\beta,\psi)$. Now the conditions characterising Landstad algebra take the form\vs

C1'. $\beta((\id{K}\tens\,\sigma)\pitilde(d))=\I{\Ahat}\tens\,(\id{K}\tens\,\sigma)\pitilde(d)$ for any $d\in D$,\vs

C2'. $B=\psi(\Ahat)(\id{K}\tens\,\sigma)\pitilde(D)$,\vs

C3'. The $\cst$-algebra generated by 
\[
\etyk{0456}
\Vset{(\id{B}\tens\,\omega)\left(V_{\psi3}(\id{K}\tens\,\sigma)\pitilde(D)_{12}V_{\psi3}^{*}\right)}{\Vs{4}\omega\in A_{*}}
\]
coincides with $(\id{K}\tens\,\sigma)\pitilde(D)$.\vs

Adapting Condition C3 to the present context we had to replace $V_{\psi2}$ by $V_{\psi3}$, because $\psi$ itself stays for two legs: $\psi(\ahat)=\I{K}\tens\sigmahat(\ahat)\in\B(K\tens H)$. Clearly $V_{\psi3}=V_{\sigmahat3}$.\vs

 Conditions C2' and C1' are already verified (cf. \rf{75} and the second formula of \rf{1632}). Formula \rf{1708} shows that the set \rf{0456} coincides with $(\id{K}\tens\,\sigma)\pitilde(D_{1})$, where $D_{1}$ is given by \rf{910}. Now Condition C3' follows directly from Condition WA1 of Definition \ref{mea}. Hence $\pitilde(D)$ is the Landstad algebra of $(B,\beta,\psi)$. 
\vs

$\cst$-algebra $D$ and the Landstad algebra of $(B,\beta,\psi)$ are obviously isomorphic: $(\id{K}\tens\,\sigma)\pitilde$ is the isomorphism. Let $((\id{K}\tens\,\sigma)\pitilde(D),\alpha')$ be the Landstad $G$-dynamical system of $(B,\beta,\psi)$. Now formula \rf{1903} takes the form
\[
\alpha'(\omega,(\id{K}\tens\,\sigma)\pitilde(d))=(\id{B}\tens\,\omega)\left(V_{\psi3}\pitilde(d)_{1\sigma}V_{\psi3}^{*}\right)
\] 
for any $\omega\in A_{*}$ and $d\in D$. Remembering that $V_{\psi3}=V_{\sigmahat3}$ and taking into account \rf{1708} we get
\[
\alpha'(\omega,(\id{K}\tens\,\sigma)\pitilde(d))=(\id{K}\tens\,\sigma)\pitilde(\alpha(\omega,d)).
\]
In other words $(\id{K}\tens\,\sigma)\pitilde$ is an isomorphism of weak $G$-dynamical systems. We see that the Landstad $G$-dynamical system of $(B,\beta,\psi)$ and the original weak $G$-dynamical system $(D,\alpha)$ are isomorphic. The proofs of statement 3 and of Theorem \ref{CP} are complete. 
\end{pf}\vs
%\newpage

\section{Weak action implemented by a unitary representation.}
\label{sek8}
Let $C$ be a $\cst$ algebra and $U$ be a unitary element of $\M(C\tens A)$. We say that $U$ is a unitary representation of $G$ if
\[
\etyk{BW1}
(\id{C}\tens\,\Delta)U=U_{12}U_{13}.
\]

Interesting examples of weak $G$-dynamical systems are provided by the following Theorem:

\begin{Thm}
\label{TO}
Let $C$ is a $\cst$-algebra, $U\in\M(C\tens A)$ be a unitary representation of $G$ and
%\[
%X=(\I{C}\tens\,A)U(C\tens\I{A})\subset\M(C\tens A).
%\]
\[
D_{1}=\Vset{(\id{C}\tens\,\omega)\left(U(C\tens\I{A})U\Vs{3}^{*}\right)}{\omega\in A_{*}}
\]
%Then $XX^{*}$ is of the form
%\[
%XX^{*}=D_{1}\tens A,
%\]
Then $D_{1}$ is a norm-closed subspace of $\M(C)$. Let $D\subset\M(C)$ be a $\cst$-subalgebra generated by $D_{1}$. For any $\omega\in A_{*}$ and $d\in D$ we set
\[
\etyk{3100}
\alpha(\omega,d)=(\id{C}\tens\,\omega)\left(U(d\tens\I{A})U^{*}\right).
\]
Then $\alpha$ is a right weak action of $G$ on $D$.
\end{Thm}

We say that $\alpha$ is a weak action implemented by $U$. The reader should notice that $D$ need not coincide with $C$. Corresponding $G$-product is described by the following Theorem:

\begin{Thm}
\label{TO1}
Let $C$ is a $\cst$-algebra and $U\in\M(C\tens A)$ be a unitary representation of $G$. Then there exists unique morphism $\psi\in\Mor(\Ahat,\Ahat\tens C)$ such that 
\[
\etyk{BW4}
(\psi\tens\id{A})V=U_{23}V_{13},
\]
It makes the diagram
\[
\etyk{BW5}
\vcenter{
\xymatrix{
\Ahat\ar[rr]^-{\psi}\ar[d]_-{\Deltahat}&&\Ahat\tens C\ar[d]^-{\Deltahat\tens\id{C}}\\
\Ahat\tens\Ahat\ar[rr]_-{\id{\Ahat}\tens\,\psi}&&\Ahat\tens\Ahat\tens C
}}
\]
commutative. Comparing this diagram with \rf{a1} we see that $(\Ahat\tens C,\Deltahat\tens\id{C},\psi)$ is a $G$-product.
\end{Thm}

With the notation introduced in the above Theorems we have

\begin{Prop}
\label{TO2}
Let $(D',\alpha')$ be the Landstad $G$-dynamical system related to $(\Ahat\tens C,\Deltahat\tens\id{C},\psi)$. Then
\[
D'=\I{\Ahat}\tens\,D
\]
and
\[
\etyk{cc}
\alpha'(\omega,\I{\Ahat}\tens\,d)=\I{\Ahat}\tens\,\alpha(\omega,d)
\]
for any $\omega\in A_{*}$ and $d\in D$.
\end{Prop}

At first we shall prove Theorem \ref{TO1} and Proposition \ref{TO2}. Then Theorem \ref{TO} will follow from an obvious isomorphism connecting $(D,\alpha)$ and $(D',\alpha')$.

\begin{pf}[ of Theorem \ref{TO1}] We chose a Heisenberg pair acting on a Hilbert space $H$. %Let $(\sigma,\sigmahat)$ be a Heisenberg pair acting on a Hilbert space $H$. 
Applying $\id{C}\tens\,\sigma\tens\id{A}$ to the both sides of \rf{BW1} and using first formula of \rf{impleme} we get
$
U_{1\sigma}U_{13}=\left(\Vs{3}\id{C}\tens\,(\sigma\tens\id{A})\Delta\right)U
=V_{\sigmahat3}U_{1\sigma}V\Vs{3}^{*}_{\sigmahat3}
$.
Therefore 
\[
U\Vs{3}^{*}_{1\sigma}V_{\sigmahat3}U_{1\sigma}=
U_{13}V_{\sigmahat3}
\]
For any $\ahat\in\Ahat$ we set
\[
%\etyk{BW2}
\psitilde(\ahat)=
\flip\left(U\Vs{3}^{*}_{1\sigma}(\I{C}\tens\,\sigmahat(\ahat))U_{1\sigma}\right).
\]
Then $U_{1\sigma},\I{C}\tens\,\sigmahat(\ahat)\in\M(C\tens\B_{0}(H))$,  $\psitilde(\ahat)\in\M(B_{0}(H)\tens C)$ and $\psitilde\in\Mor(\Ahat,B_{0}(H)\tens C)$. We have
\[
\etyk{3605}
(\psitilde\tens\id{A})V=
\flip_{12}\left(U\Vs{3}^{*}_{1\sigma}V_{\sigmahat3}U_{1\sigma}\right)
=\flip_{12}\left(\Vs{3}U_{13}V_{\sigmahat3}\right)=U_{23}V_{\sigmahat3}.
\]
The reader should notice that in the above computation $U_{13}V_{\sigmahat3}=(\id{C}\tens\,\sigmahat\tens\id{A})U_{13}V_{23}$, whereas $U_{23}V_{\sigmahat3}=(\sigmahat\tens\id{C}\tens\id{A})U_{23}V_{13}$. Clearly $U_{23}V_{\sigmahat3}\in\M(\B_{0}(H)\tens C\tens A)$. More precisely $U_{23}V_{\sigmahat3}$ belongs to $\M(\sigmahat(\Ahat)\tens C\tens A)$. Now Theorem 1.6.6 of \cite{SLW96c} shows that $\psitilde\in\Mor(\Ahat,\sigmahat(\Ahat)\tens C)$. We know that $\sigmahat$ is faithful. Therefore $\psitilde$ is of the form $\psitilde=(\sigmahat\tens\id{C})\psi$, where $\psi\in\Mor(\Ahat,\Ahat\tens C)$. With this notation \rf{3605} implies \rf{BW4}.\vs

Taking into account \rf{BW4} we obtain 
\[
\begin{array}{r@{\;=\;}l}
\left((\Deltahat\tens\id{C})\psi\tens\id{A}\right)V&
\left(\Deltahat\tens\id{C}\tens\id{A}\right)U_{23}V_{13}=U_{34}V_{24}V_{14}\\
\Vs{5}&
\left(\id{\Ahat}\tens\,\psi\tens\id{A}\right)V_{23}V_{13}=
\left((\id{\Ahat}\tens\,\psi)\Deltahat\tens\id{A}\right)V.
\end{array}
\]
Now second formula of \rf{algabry} shows that \rf{BW5} is a commutative diagram.
Using coassociativity of $\Deltahat$ and the cancelation property one can easily show that $\Deltahat\tens\id{C}\in\Mor(\Ahat\tens C,\Ahat\tens\Ahat\tens C)$ is a continuous left action of $\Ghat$ on $\Ahat\tens C$. Hence $(\Ahat\tens C,\Deltahat\tens\id{C},\psi)$ is a $G$-product.
\end{pf}\vs

\begin{pf}[ of Proposition \ref{TO2}]
Let $B=\Ahat\tens C$, $\beta=\Deltahat\tens\id{C}\in\Mor(B,\Ahat\tens B)$ and $\psi\in\Mor(\Ahat,B)$ be the morphism introduced in Theorem \ref{TO1}. We shall prove that
\[
\etyk{BW3}
B=(\I{\Ahat}\tens\,C)\psi(\Ahat).
\]
Indeed
\[
\psi(\Ahat)=\psi\left(\Vset{(\id{\Ahat}\tens\,\omega)V\Vs{3.5}}{\omega\in A_{*}}\right)=
\Vset{(\id{\Ahat\tens C}\tens\,\omega)
U_{23}V_{13}\Vs{3.5}}
{\omega\in A_{*}}
\]
Therefore
\[
\begin{array}{r@{\;=\;}l}
(\I{\Ahat}\tens C)\psi(\Ahat)&
\Vset{(\id{\Ahat}\tens\id{C}\tens\,\omega)
\left((\I{\Ahat}\tens C\tens\I{A})U_{23}V_{13}\Vs{3.5}\right)}
{\omega\in A_{*}}\\
\Vs{7}&
\Vset{(\id{\Ahat}\tens\id{C}\tens\,\omega)
\left((\I{\Ahat}\tens C\tens A)U_{23}V_{13}\Vs{3.5}\right)}
{\omega\in A_{*}}\\
\Vs{7}&
\Vset{(\id{\Ahat}\tens\id{C}\tens\,\omega)
\left((\I{\Ahat}\tens C\tens A)V_{13}\Vs{3.5}\right)}
{\omega\in A_{*}}\\
\Vs{7}&
\Vset{(\id{\Ahat}\tens\,\omega)V\Vs{3.5}}
{\omega\in A_{*}}\tens C=\Ahat\tens C
\end{array}
\]
and \rf{BW3} follows. In the above computation $\omega$ emits $A$ to the right, next $U$ is absorbed by $C\tens A$ (this is because $U$ is a unitary element of $\M(C\tens A)$) and finally $\omega$ absorbs $A$. 
%Using \rf{BW3} we get: $B\psi(\Ahat)=(\I{\Ahat}\tens C)\psi(\Ahat)\psi(\Ahat)=(\I{\Ahat}\tens C)\psi(\Ahat)=B$. It shows that $\psi\in\Mor(\Ahat,B)$. 
\vs 

To determine the Landstad algebra of $(\Ahat\tens C,\Deltahat\tens\id{C},\psi)$ we shall use the procedure described in Section \ref{sek5}. Let $(\sigma,\sigmahat)$ be a Heisenberg pair acting one Hilbert space $H$ and $\varphi\in\Mor(B,B\tens\B_{0}(H))$ be the morphism introduced by \rf{859}. We have
\[
\begin{array}{r@{\;=\;}l}
\beta(\I{\Ahat}\tens\,C)&\I{\Ahat}\tens\I{\Ahat}\tens\,C,\\
\beta'(\I{\Ahat}\tens\,C)\Vs{5}&\I{\Ahat}\tens\,C\tens\I{\Ahat},\\
\varphi(\I{\Ahat}\tens\,C)\Vs{5}&V_{\psi\sigma}(\I{\Ahat}\tens\,C\tens\I{H})V\Vs{3}_{\psi\sigma}^{*}.
\end{array}
\]
Formula \rf{BW4} says that $V_{\psi3}=U_{23}V_{13}$. Hence $V_{\psi\sigma}=U_{2\sigma}V_{1\sigma}$. Second leg of $V_{1\sigma}$ is trivial. Therefore $V_{1\sigma}$ commutes with $\I{\Ahat}\tens\,C\tens\I{H}$ and
\[
\begin{array}{r@{\;=\;}l}
\varphi(\I{\Ahat}\tens\,C)&U_{2\sigma}(\I{\Ahat}\tens\,C\tens\I{H})U\Vs{3}_{2\sigma}^{*}\\ \Vs{5}&\I{\Ahat}\tens\,(\id{C}\tens\,\sigma)\left(U(C\tens\I{A})U\Vs{3}^{*}\right).
\end{array}
\]
Formula \rf{1731f} shows that $\varphi(\psi(\Ahat))=\I{B}\tens\,\sigmahat(\Ahat)$. Taking into account \rf{BW3} we obtain
\[
\varphi(B)=\I{\Ahat}\tens\,(\id{C}\tens\,\sigma)\left(U(C\tens\I{A})U\Vs{3}^{*}\right)
\left(\I{C}\tens\,\sigmahat(\Ahat)\right).
\]
In this section $D_{1}$ and $D$ denote operator spaces introduced in Theorem \ref{TO}. To avoid conflict of notation with the one used earlier we decorate $D_{1}$ in \rf{937} with prime:
\[
\begin{array}{r@{\;=\;}l}
D'_{1}&\Vset{(\id{B}\tens\,\mu)\varphi(B)}{\Vs{3.5}\mu\in\B(H)_{*}}
\\
\Vs{6}&\I{\Ahat}\tens\Vset{(\id{C}\tens\,\mu)
\left((\id{C}\tens\,\sigma)\left(U(C\tens\I{A})U\Vs{3}^{*}\right)\left(\I{C}\tens\,\sigmahat(\Ahat)\right)\right)}{\mu\in\B(H)_{*}}.
\end{array}
\]
Absorbing $\sigmahat(\Ahat)$ by $\mu$ and replacing $\mu\comp\sigma$ by $\omega$ we obtain
\[
D'_{1}=\I{\Ahat}\tens\Vset{(\id{C}\tens\,\omega)
\left(U(C\tens\I{A})U\Vs{3}^{*}\right)}{\omega\in A_{*}}=\I{\Ahat}\tens D_{1}.
\]
By Proposition \ref{LAEU}, the Landstad algebra of $(\Ahat\tens\Ahat,\Deltahat\tens\id{C},\psi)$ coincides with $\I{\Ahat}\tens D$. According to \rf{BW4} $V_{\psi3}=U_{23}V_{13}$. Now, formula \rf{1903} takes the form
\[
\begin{array}{r@{\;=\;}l}
\alpha'(\omega, \I{\Ahat}\tens\,d)&(\id{\Ahat}\tens\id{C}\tens\,\omega)\left(U_{23}V_{13}(\I{\Ahat}\tens\,d\tens\I{A})V_{13}^{*}U_{23}^{*}\right)\\
\Vs{6}&\I{\Ahat}\tens\,(\id{C}\tens\,\omega)\left(U(d\tens\I{A})U^{*}\right)=\I{\Ahat}\tens\,\alpha(\omega,d)
\end{array}
\]
for any $\omega\in A_{*}$ and $d\in D$.
\end{pf}\vs

\begin{pf}[ of Theorem \ref{TO}]
Let $\imath(d)=\I{\Ahat}\tens\,d$ any for $d\in D$. Then $\imath:D\longrightarrow\I{\Ahat}\tens\,D=D'$ is an isomorphism of $\cst$-algebras. Formula \rf{cc} takes the form 
\[
\alpha'(\omega,\imath(d))=\imath(\alpha(\omega,d)).
\]
It shows that $(D',\alpha')$ and $(D,\alpha)$ are isomorphic. By Proposition \rf{TO2}, $(D',\alpha')$ is a weak $G$-dynamical system. So is $(D,\alpha)$.
\end{pf}\vs

We say that $U$ is regular if $(\I{C}\tens A)U(C\tens\I{A})=C\tens A$. If $U$ is regular then $D=D_{1}=C$. In general $D\neq C$. To obtain the most obvious example of the above construction we set: $C=\Ahat$ and $U=V$. Then $\psi=\Deltahat$ and diagram \rf{BW5} takes the form
\[
\xymatrix{
\Ahat\ar[rr]^-{\Deltahat}\ar[d]_-{\Deltahat}&&\Ahat\tens\Ahat\ar[d]^-{\Deltahat\tens\id{\Ahat}}\\
\Ahat\tens\Ahat\ar[rr]_-{\id{\Ahat}\tens\,\Deltahat}&&\Ahat\tens\Ahat\tens\Ahat
}
\]
This is the commutative diagram stating the coassociativity of $\Deltahat$. This way for any locally compact quantum group $G$ we have canonically associated $G$-product $(\Ahat\tens\Ahat,\Deltahat\tens\id{\Ahat},\Deltahat)$. Now
\[
\etyk{ost}
D_{1}=\Vset{(\id{\Ahat}\tens\,\omega)\left(V(\Ahat\tens\I{H})V^{*}\right)}{\omega\in\B(H)_{*}}
\]
and $D$ is the $\cst$-algebra generated by $D_{1}$. $D_{1}$ and $D$are subsets of $\M(\Ahat)$. For any $\omega\in A_{*}$ and $d\in D$ we have
\[
\alpha(\omega,d)=(\id{\Ahat}\tens\,\omega)\left(V(d\tens\I{H})V^{*}\right).
\]
$\alpha$ is a right weak action of $G$ on $D$ and $(D,\alpha)$ is a $G$-dynamical system corresponding to the $G$-product $(\Ahat\tens\Ahat,\Deltahat\tens\id{\Ahat},\Deltahat)$. If $G$ is regular then $D=\Ahat$. Otherwise $D'\neq\Ahat$.\vs

The first known example of non-regular locally compact quantum group was quantum deformation $E_{q}(2)$ of the group of motions of Euclidean plane (with real deformation parameter $0<q<1$, see \cite{SLW91a, SLW91d} and \cite{BaajE2}). Let\footnote{notice change of notation: the role of $G$ and $\Ghat$ are interchanged.} $A={\rm C}_{0}(E_{q}(2))$ and $\Delta\in\Mor(A,A\tens A)$ be the corresponding comultiplication: $E_{q}(2)=(A,\Delta)$. Then $(A\tens A,\Delta\tens\id{A},\Delta)$ is a $\widehat{E_{q}(2)}$-product. Let $D'=\I{A}\tens D$ be its Landstad algebra. Condition C2 says that 
\[
\Delta(A)(\I{A}\tens D)=A\tens A.
\] 
A boring calculations entering deeply into anatomy of $E_{q}(2)$ shows that $D$ is unital. Therefore $\Delta(A)\subset A\tens A$. For the first time this unexpected result appeared in \cite{SLW92c}.%\newpage
 
\section{Kasprzak approach to Rieffel deformation.}
\label{sek6}

The following example comes from the Kasprzak theory \cite{PK09, PK10, PK11}. He was able to find an elegant realisation of Rieffel deformation \cite{Rief1} of $\cst$-algebras endowed with an action of a group. Kasprzak (and Rieffel) worked with locally compact abelian group, but due to the further developement we may consider any locally compact quantum group $G$. We shall use the notation introduced in previous sections. \vs

Assume that we have a unitary two-cocycle. This is a unitary element $\Omega\in\M(\Ahat\tens\Ahat)$ such that
\[
\etyk{2co}
(\Omega\tens\I{\Ahat})(\Deltahat\tens\id{\Ahat})(\Omega)=(\I{\Ahat}\tens\Omega)(\id{\Ahat}\tens\,\Deltahat)(\Omega).
\]
We shall also assume that the Drinfeld twist induced by $\Omega$ is trivial:
\[
\etyk{DT}
\Omega^{*}\Deltahat(\ahat)\Omega=\Deltahat(\ahat)
\]
for any $\ahat\in\Ahat$.\vs

Let $D$ be a $\cst$-algebra equipped with a right weak action $\alpha$ of $G$. Then $(D,\alpha)$ be a $G$-dynamical system. Using  Theorem \ref{CP} we may find $G$-product $(B,\beta,\psi)$ with Landstad dynamical system isomorphic to $(D,\alpha)$.

\begin{Thm}
\label{Kasp}
For any $b\in B$ we set
\[
\beta^{\Omega}(b)=\Omega_{1\psi}^{*}\beta(b)\Omega_{1\psi}. 
\]
Then $\beta^{\Omega}(b)\in\M(\Ahat\tens B)$, $\beta^{\Omega}\in\Mor(B,\Ahat\tens B)$ is a continuous left action of $\Ghat$ on $B$ and $(B,\beta^{\Omega},\psi)$ is a $G$-product.
\end{Thm} 
\begin{pf}
We have to show that
\begin{align}
(\id{\Ahat}\tens\beta^{\Omega})\beta^{\Omega}(b)&=(\Deltahat\tens\id{B})\beta^{\Omega}(b), \label{st}\\
(\Ahat\tens\I{B})\beta^{\Omega}(B)&=\Ahat\tens B, \label{stst}\\
\beta^{\Omega}(\psi(\ahat))&=(\id{\Ahat}\tens\,\psi)\Deltahat(\ahat). \label{ststst}
\end{align}

Relation \rf{ststst} is easy to verify:
\[
\begin{array}{r@{\;=\;}l}
\beta^{\Omega}(\psi(\ahat))&\Omega_{1\psi}^{*}\beta(\psi(\ahat))\Omega_{1\psi}
=\Omega_{1\psi}^{*}(\id{\Ahat}\tens\,\psi)\Deltahat(\ahat))\Omega_{1\psi}\\
\Vs{5}&(\id{\Ahat}\tens\,\psi)\left(\Omega^{*}\Deltahat(\ahat))\Omega\right)
=(\id{\Ahat}\tens\,\psi)\Deltahat(\ahat).
\end{array}
\]

We already know that $\Omega_{1,\beta^{\Omega}\psi}=(\id{\Ahat}\tens\beta^{\Omega}\psi)\Omega=
(\id{\Ahat}\tens\,(\id{\Ahat}\tens\,\psi)\Deltahat)\Omega=\Omega_{1,(\id{\Ahat}\tens\,\psi)\Deltahat}$. We com\-pute
\[
(\id{\Ahat}\tens\beta^{\Omega})\beta^{\Omega}(b)=\Omega_{1,\beta^{\Omega}\psi}^{*}(\id{\Ahat}\tens\beta^{\Omega})\beta(b)\Omega_{1,\beta^{\Omega}\psi}=
\Omega_{1,\beta^{\Omega}\psi}^{*}\Omega_{2\psi}^{*}(\id{\Ahat}\tens\beta)\beta(b)\Omega_{2\psi}\Omega_{1,\beta^{\Omega}\psi}.
\]
On the other hand
\[
(\Deltahat\tens\id{B})\beta^{\Omega}(b)=\Omega_{\Deltahat,\psi}^{*}(\Deltahat\tens\id{B})\beta(b)\Omega_{\Deltahat,\psi}.
\]
To prove \rf{st} it is enough to show that $\Omega_{2\psi}\Omega_{1,\beta^{\Omega}\psi}\Omega_{\Deltahat,\psi}^{*}$ commutes with $(\Deltahat\tens\id{B})\beta(b)$. Using at the last moment \rf{2co} we obtain
\[
\begin{array}{r@{\;=\;}l}
\Omega_{2\psi}\Omega_{1,\beta^{\Omega}\psi}\Omega_{\Deltahat,\psi}^{*}&
\Omega_{2\psi}\Omega_{1,(\id{\Ahat}\tens\,\psi)\Deltahat}\Omega_{\Deltahat,\psi}^{*}\\
\Vs{5}&(\id{\Ahat}\tens\id{\Ahat}\tens\,\psi)\left(\Vs{4}
(\I{\Ahat}\tens\Omega)(\id{\Ahat}\tens\,\Deltahat)\Omega(\Deltahat\tens\id{\Ahat})\Omega^{*}
\right)\\
\Vs{5}&\Omega\tens\I{B}.
\end{array}
\]
Now the commutativity follows immediately form \rf{DT}.\vs

We shall show \rf{stst}. $\Omega$ is a unitary multiplier of $\Ahat\tens\Ahat$. Therefore $(\Ahat\tens\Ahat)\Omega^{*}=\Ahat\tens\Ahat$. Using the cancelation formula $\Ahat\tens\Ahat=(\Ahat\tens\I{\Ahat})\Deltahat(\Ahat)$ and relation \rf{DT} we obtain
\[
(\Ahat\tens\I{\Ahat})\Omega^{*}\Deltahat(\Ahat)=(\Ahat\tens\I{\Ahat})\Deltahat(\Ahat).
\]
Applying to the both sides $\id{\Ahat}\tens\,\psi$ and taking into account commutativity of \rf{a1} we get
\[
(\Ahat\tens\I{B})\Omega_{1\psi}^{*}\beta(\psi(\Ahat))=
(\Ahat\tens\I{B})\beta(\psi(\Ahat)).
\]
We know that $\psi\in\Mor(\Ahat,B)$. Therefore $\psi(\Ahat)B=B$. Multiplying both sides of the above formula by $\beta(B)$  and using Podleś condition for the action $\beta$ we have
\[
(\Ahat\tens\I{B})\Omega_{1\psi}^{*}\beta(B)=
(\Ahat\tens\I{B})\beta(B)=\Ahat\tens B.
\]
Finally
\[
(\Ahat\tens\I{B})\beta^{\Omega}(B)=(\Ahat\tens B)\Omega_{1\psi}=\Ahat\tens B.
\]
\end{pf}\vs

Let $D'$ be the Landstad algebra related to $G$-product $(B,\beta^{\Omega},\psi)$. According to Kasprzak, $D'$ may be considered as Rieffel deformation of $D$. Recently Kasprzak theory was extended by Neshveyev and Tuset \cite{NeshTuset} for non-trivial Drinfeld twist i.e. when \rf{DT} does not hold. Then they had to consider deformations of $G$ and $\psi$.

%\newpage

\section{A remark on Landstad conditions}
\label{sek7}

This section is not in the main stream of the paper. It contains a generalisation of a result of Kasprzak that was used to simplify second Landstad condition (cf. formulae (3) and (4) of \cite{PK09}). In this section we assume that $A$ admits a continuous counit. (i.e: $G$ is coameanable).

\begin{Thm}
\label{0925}
Let $(B,\beta,\psi)$ be a $G$-product and $d\in\M(B)$. Assume that
\[
\etyk{0926}
V_{\psi2}(d\tens\I{A})V_{\psi2}^{*}(\I{B}\tens a)\in\M(B)\tens A
\]
for any $a\in A$. Then the following conditions are equivalent:
\begin{enumerate}
\item $\psi(\ahat)d\in B$ for any $\ahat\in\Ahat$,
\item $d\psi(\ahat')\in B$ for any $\ahat'\in\Ahat$,
\item $\psi(\ahat)d\psi(\ahat')\in B$ for any $\ahat,\ahat'\in\Ahat$.
\end{enumerate}
\end{Thm}
We shall use the following result:
\begin{Prop}
\label{0972}
There exists a bounded net $(\ehat_{\lambda})_{\lambda\in\Lambda}$ of elements of $\Ahat$ converging strictly to $\I{\Ahat}$ such that
\[
\etyk{1125}
\lim_{\lambda\in\Lambda}\norm{\left[d\Vs{3.2},\psi(\ehat_{\lambda})\right]}=0
\]
for any $d\in\M(B)$ satisfying relation \rf{0926}. Square bracket in the above formula denotes commutator.
\end{Prop}
\begin{pf}
Let $e$ be a counit of $A$ and $(\omega_{\lambda})_{\lambda\in\Lambda}$ be a net of normal states on $A$ weakly converging to $e$: 
\[
\lim_{\lambda\in\Lambda}\omega_{\lambda}(x)=e(x)
\]
for any $x\in A$. Then for any $r\in\M(B)\tens A$ we have
\[
\etyk{0256}
\text{norm-}\!\lim_{\lambda\in\Lambda}(\id{B}\tens\,\omega_{\lambda})r=(\id{B}\tens\, e)r.
\]
Indeed \rf{0256} is obvious for $r\in\M(B)\tens\hspace{-.2em}_{\text{alg}}\, A$. Moreover the set of all $r\in\M(B)\tens A$ satisfying  \rf{0256} is closed in norm topology. Remembering that $\M(B)\tens A$ is a norm closure of $\M(B)\tens\hspace{-.2em}_{\text{alg}}\, A$ we obtain \rf{0256} in full generality.\vs

We fix an element $a\in A$ such that $e(a)=1$ and set
\[
\etyk{0325}
\ehat_{\lambda}=(\id{B}\tens\,\omega_{\lambda})\left(V^{*}(\I{\Ahat}\tens\,a)\right).
\] 
Then $\ehat_{\lambda}\in\Ahat$ for any $\lambda\in\Lambda$. With an easy calculation we obtain
\[
\left[d\Vs{3.2},\psi(\ehat_{\lambda})\right]=(\id{B}\tens\,\omega_{\lambda})X,
\]
where
\[
\begin{array}{r@{\;=\;}l}
X&(d\tens\I{A})V_{\psi2}^{*}(\I{B}\tens\,a)-V_{\psi2}^{*}(d\tens a)\\
\Vs{6}&V_{\psi2}^{*} \left[\Vs{3.2} V_{\psi2}(d\tens\I{A})V_{\psi2}^{*}(\I{B}\tens a)-d\tens a\right].
\end{array}
\]
Assume that $d\in\M(B)$ satisfies \rf{0926}. Then the expression in square bracket belongs to $\M(B)\tens A$. Therefore $X^{*}X\in\M(B)\tens A$. It is known that $(\id{B}\tens\,e)V=\I{\Ahat}$. Therefore ($e$ is a character) $(\id{B}\tens\,e)X=d-d=0$ and $(\id{B}\tens\,e)(X^{*}X)=0$. Formula \rf{0256} shows now that
\[
\lim_{\lambda\in\Lambda}\norm{(\id{B}\tens\,\omega_{\lambda})(X^{*}X)}=0.
\]
Clearly $\id{B}\tens\,\omega_{\lambda}$ is a completely positive unital mapping. Using Kadison inequality (cf \cite[Corollary 1.3.2, page 9]{Stormer13}) we obtain
\[
\begin{array}{rl}
\norm{\left[d\Vs{3.2},\psi(\ehat_{\lambda})\right]}^{2}&=\norm{(\id{B}\tens\,\omega_{\lambda})(X^{*})(\id{B}\tens\,\omega_{\lambda})(X)}\\
\Vs{5}&\leq\norm{(\id{B}\tens\,\omega_{\lambda})(X^{*}X)}.
\end{array}
\]
Formula \rf{1125} is shown. To end the proof we have to show that the net $(\ehat_{\lambda})_{\lambda\in\Lambda}$ converges strictly to $\I{\Ahat}\in\M(\Ahat)$. To this end we have to choose the net $(\omega_{\lambda})_{\lambda\in\Lambda}$ in a more specific way.\vs

Let $(\omega'_{\lambda})_{\lambda\in\Lambda}$ be a net of normal states on $A$ weakly converging to $e$ and $c$ be an element of $A$ such that $e(c)=1$. For any $a\in A$ and $\lambda\in\Lambda$ we set
\[
\omega_{\lambda}(a)=\frac{\omega'_{\lambda}(c^{*}ac)}{\omega'_{\lambda}(c^{*}c)}.
\]
Then $(\omega_{\lambda})_{\lambda\in\Lambda}$ is a net of normal states on $A$ weakly converging to $e$. Now \rf{0325} takes the form
\[
\ehat_{\lambda}=\frac{(\id{B}\tens\,\omega'_{\lambda})\left((\I{\Ahat}\tens c^{*})V^{*}(\I{\Ahat}\tens ca)\right)}{\omega'_{\lambda}(c^{*}c)}
\] 
and
\[
\begin{array}{r@{\;=\;}l}
\ahat\ehat_{\lambda}&{\displaystyle \frac{(\id{B}\tens\,\omega'_{\lambda})\left((\ahat\tens c^{*})V^{*}(\I{\Ahat}\tens ca)\right)}{\omega'_{\lambda}(c^{*}c)}},\\
\Vs{8}\ehat_{\lambda}\ahat&{\displaystyle \frac{(\id{B}\tens\,\omega'_{\lambda})\left((\I{\Ahat}\tens c^{*})V^{*}(\ahat\tens ca)\right)}{\omega'_{\lambda}(c^{*}c)}}
\end{array}
\] 
for any $\ahat\in\Ahat$. The reader should notice that $(\ahat\tens c^{*})V^{*}(\I{\Ahat}\tens ca)$ and $(\I{\Ahat}\tens c^{*})V^{*}(\ahat\tens ca)$ belong to $B\tens A$. Formula \rf{0256} shows now that $\ahat\ehat_{\lambda}$ and $\ehat_{\lambda}\ahat$ converge in norm to $\ahat$. It means that $\ehat_{\lambda}$ converges strictly to $\I{\Ahat}\in\M(\Ahat)$.
\end{pf}

\begin{pf}[ of Theorem \ref{0925}] 
We know that $\psi(\ahat),\psi(\ahat')\in\M(B)$. Therefore (3) follows from (1) and from (2). We shall prove that (1) follows from (3). One can easily verify that 
\[
\psi(\ahat)d-\psi(\ahat)d\psi(\ehat_{\lambda})=\psi(\ahat-\ahat\ehat_{\lambda})d-\psi(\ahat)\left[d\Vs{3.2},\psi(\ehat_{\lambda})\right].
\]
Proposition \ref{0972} shows now that
\[
\psi(\ahat)d=\text{norm-}\!\lim_{\lambda\in\Lambda}\psi(\ahat)d\psi(\ehat_{\lambda}).
\]
Remembering that $B$ is norm-closed in $\M(B)$ we see that (3) implies (1). In the similar way one shows that (3) implies (2) 
\end{pf}

\section*{Acknowlegdement}
The research started in Oberwolfach within the program "research in pairs" in August 2014 and ended in Warsaw in last week of June 2017. The authors are very grateful to the Oberwolfach Research Institute for Mathematics and to Banach Center in Warsaw for creating the perfect conditions for fruitful work.

\vspace{1cm}

\bibliographystyle{plain}
%\bibliography{SLWbibMR,Spis}
\def\cst{{\textup{C}^*}} \def\cprime{$'$}
  \def\polhk#1{\setbox0=\hbox{#1}{\ooalign{\hidewidth
  \lower1.5ex\hbox{`}\hidewidth\crcr\unhbox0}}} \def\cprime{$'$}
  \def\cst{{\textup{C}^*}}

\vspace{5mm}

\end{document}